\documentclass[11pt]{article}

\usepackage{grffile}

\usepackage{amsmath,amsfonts,amscd,a4wide}
\usepackage{amsmath}
\usepackage{amssymb}
\usepackage{graphicx}
\usepackage{amsthm}
\usepackage{color}
\usepackage{setspace} 

\usepackage{transparent}
\usepackage[small]{caption}
\usepackage{subcaption}
\usepackage{enumitem}

\usepackage{color,soul}
\definecolor{lightgray}{gray}{0.85}
\sethlcolor{lightgray}
\usepackage[margin=1in]{geometry}

\newtheorem{Lemma}{Lemma}[section]
\newtheorem{Theorem}{Theorem}
\newtheorem{Proposition}[Lemma]{Proposition}

\newtheorem{Remark}[Lemma]{Remark}
\newtheorem{Definition}[Lemma]{Definition}

\newenvironment{Proof}%
 {\begin{trivlist} \item[]{\bf Proof. }}%
 {\hspace*{\fill}$\rule{.4\baselineskip}{.4\baselineskip}$\end{trivlist}}

\newenvironment{Acknowledgment}%
 {\begin{trivlist}\item[]\textbf{Acknowledgments.}}{\end{trivlist}}

\makeatletter\@addtoreset{figure}{section}\makeatother

\makeatletter \@addtoreset{equation}{section} \makeatother

\newcommand{\R}{\mathbb{R}}
\newcommand{\C}{\mathbb{C}}

\newcommand{\Z}{\mathbb{Z}}

	\newcommand{\mc}[1]{\mathcal{#1}}
	\newcommand{\mb}[1]{\mathbb{#1}}
	\newcommand{\tl}[1]{\tilde{#1}}
	\newcommand{\lp}{\left}
	\newcommand{\rp}{\right}
	\newcommand{\beq}{\begin{equation}}
	\newcommand{\eeq}{\end{equation}}
	\newcommand{\ba}{\begin{align}}
	\newcommand{\ea}{\end{align}}
	\newcommand{\fr}[2]{\frac{#1}{#2}}
	\newcommand{\p}{\partial}
	\newcommand{\ri}{\mathrm{i}}
	\newcommand{\rmi}{\mathrm{i}}
	\newcommand{\rlin}{\mathrm{lin}}
	\newcommand{\rmd}{\mathrm{d}}
	\newcommand{\re}{\mathrm{e}}
	\newcommand{\rme}{\mathrm{e}}
	\newcommand{\rmo}{{\scriptstyle\mathcal{O}}}
	\newcommand{\rmO}{\mathcal{O}}
	\renewcommand{\Re}{\mathrm{Re}}
	\renewcommand{\Im}{\mathrm{Im}}
	\newcommand{\rre}{\mathrm{Re}}
	\newcommand{\rim}{\mathrm{Im}}
	\newcommand{\rab}{\mathrm{abs}}
	
		\newcommand{\rtf}{\mathrm{tf}}
				\newcommand{\rpp}{\mathrm{p}}
					\newcommand{\rb}{\mathrm{b}}
			\newcommand{\rs}{\mathrm{s}}
				\newcommand{\rss}{\mathrm{ss}}
					\newcommand{\rcu}{\mathrm{cu}}

\title{Triggered Fronts in the Complex Ginzburg Landau Equation}
\author{Ryan Goh and Arnd Scheel}
\begin{document}

\maketitle
\begin{abstract}
We study patterns that arise in the wake of an externally triggered, spatially propagating instability in the complex Ginzburg-Landau equation. We model the trigger by a spatial inhomogeneity moving with constant speed. In the comoving frame, the trivial state is unstable to the left of the trigger and stable to the right. At the trigger location, spatio-temporally periodic wavetrains nucleate. Our results show existence of coherent, ``heteroclinic'' profiles when the speed of the trigger is slightly below the speed of a free front in the unstable medium. Our results also give expansions for the wavenumber of wavetrains selected by these coherent fronts. A numerical comparison yields very good agreement with observations, even for moderate trigger speeds. Technically, our results provide a heteroclinic bifurcation study involving an equilibrium with an algebraically double pair of complex eigenvalues. We use geometric desingularization and invariant foliations to describe the unfolding. Leading order terms 
are determined by a condition of oscillations in a projectivized flow, which can be found by intersecting absolute spectra with the imaginary axis. 
\end{abstract}

\section{Introduction and main results}

The complex Ginzburg-Landau equation (CGL),
\beq
A_t=(1+ \ri \alpha)\Delta A + A - (1+ \ri \gamma)A|A|^2,\label{eq:cgl}
\eeq
is a prototypical model for the spontaneous emergence of self-organized, regular spatio-temporal patterns in spatially extended systems. It can be derived and justified in a universal fashion as an approximation near the onset of an oscillatory instability in dissipative systems; see for instance \cite{Aranson02,Miecgl} and references therein. As is typical for pattern-forming systems, CGL supports a variety of coherent and complex patterns. Even in parameter regimes $\alpha\sim\gamma$, when the equation is close to a gradient-like flow, there exist continua of stable periodic patterns, $A\sim \re^{ \ri (\omega(k) t-k\cdot x)}$, and coherent defects, most prominently Nozaki-Bekki holes in one space dimension and spiral waves in two space-dimension, both having vanishing amplitude $A=0$ at $x=0$. Starting with random initial data typically leads to spatial patches of coherent patterns with boundaries that slowly evolve. On the other hand, starting with spatially localized initial conditions, one observes a 
spatial invasion process that leaves distinct wavetrains in its wake, whose wavenumber does not depend on the initial condition but only on system parameters; see \cite{deelanger,deelanger1,vanSaarloos03}. While such a regular pattern selection mechanism might be desirable, either from an engineering or a natural selection perspective, it is typically difficult to achieve since it requires the preparation of a perfect unstable state and hence uniform suppression of fluctuations almost everywhere. 

\paragraph{Spatially triggered instabilities.}

Rather than quenching the system instantaneously into an unstable state at \emph{all} locations $x$, one often observes (or engineers) a progressively triggered instability, where a trivial state is unstable only in the wake of a progressing front. In the case of CGL in one space dimension, such a process can be modeled via an inhomogeneous linear driving term,
\begin{equation}\label{e:cglt}
A_t=(1+ \ri \alpha)A_{xx} + \chi A - (1+ \ri \gamma)A|A|^2,
\end{equation}
where $\chi=\chi(x-ct)$, $c>0$, and, setting $\xi = x - ct$, $\chi(\xi)\to\chi_\pm$ for $x\to\pm\infty$, with $\chi_->0>\chi_+$.  In fact, such spatially triggered instabilities arise in a variety of engineering applications and can be viewed as caricature models of patterned growth processes. In the following, we mention a few examples that have been studied in the literature.

Spontaneous nanostructures emerge on surfaces in ion beam milling due to a secondary sputtering instability \cite{Bradley88,Bradley10}. The patterning process is usually irregular due to defect nucleation and intrinsic surface instabilities. When the surface is exposed to an ion beam through a moving mask, regular, defect-free nano-scale patterns can be created \cite{bradleygelfand}. Controlling the speed and angle of the beam, a rich variety of  regular, defect free, nano-scale structures can be obtained. Similar models and effects arise in abrasive water jet cutting \cite{friedrich}, and in phase separation processes \cite{Foard12,krekhov}. In a similar vein, recurrent precipitation patterns arise in the wake of a diffusion front to form banded ring-like patterns, as originally observed by Liesegang \cite{liese}. Models often include a simple bimolecular reaction whose product feeds as a diffusive source term into a model for the precipitation kinetics; see for instance \cite{Droz00,Keller81}. Lastly, we 
mention patterns produced through growth and chemotaxis in bacterial colonies; see for instance \cite{matsuyama}. In a first approximation, one can envision patterns formed by chemotactic motion in the wake of a spatially spreading growth process \cite{mimura}, with phenomena reminiscent of patterns in the wake of trigger fronts in 2d Cahn-Hilliard equations \cite{Foard12}. 

In the present work, we focus on the complex Ginzburg-Landau equation. We believe that many of the features of our analysis can be immediately transfered to more general pattern-forming systems, such as the examples listed above. The advantage of the Ginzburg-Landau equation is that, due to the gauge symmetry, periodic patterns are explicit and invasion as well as triggered fronts can be found as heteroclinic orbits in a 4-dimensional traveling-wave equation, thus avoiding the (well-understood) complications of infinite-dimensional, ill-posed systems for modulated traveling waves \cite{Goh11,imd,kirchgassner,SansSch01,Scheel03}. On the other hand, the complex Ginzburg-Landau equation is of interest since it approximately describes many pattern-forming systems near onset.

\paragraph{Phenomenology.}

Going back to the Ginzburg-Landau eqation with a moving triggered instability (\ref{e:cglt}), we now describe the basic scenario that we will analyze quantitatively in this article. In the following work we fix
\beq
\chi(\xi) = \begin{cases}
\chi_-=+1, &\mbox{for }\xi\leq0 \\
\chi_+=-1, &\mbox{for }\xi>0 
\end{cases}\label{eq:trigger}
\eeq
and consider (\ref{e:cglt}) with fixed $\alpha,\gamma$, while varying the speed $c$ of the trigger $\chi$ as our primary bifurcation parameter. Our results can be eaily adapted to different values of $\chi_\pm$. Smooth triggers $\chi$ with derivative $\chi'$ sufficiently localized would also be immediately amenable to the following analysis. 

\begin{figure}[h!]
        \centering
        \begin{subfigure}[h!]{0.45\textwidth}
                \centering
                \includegraphics[width=\textwidth]{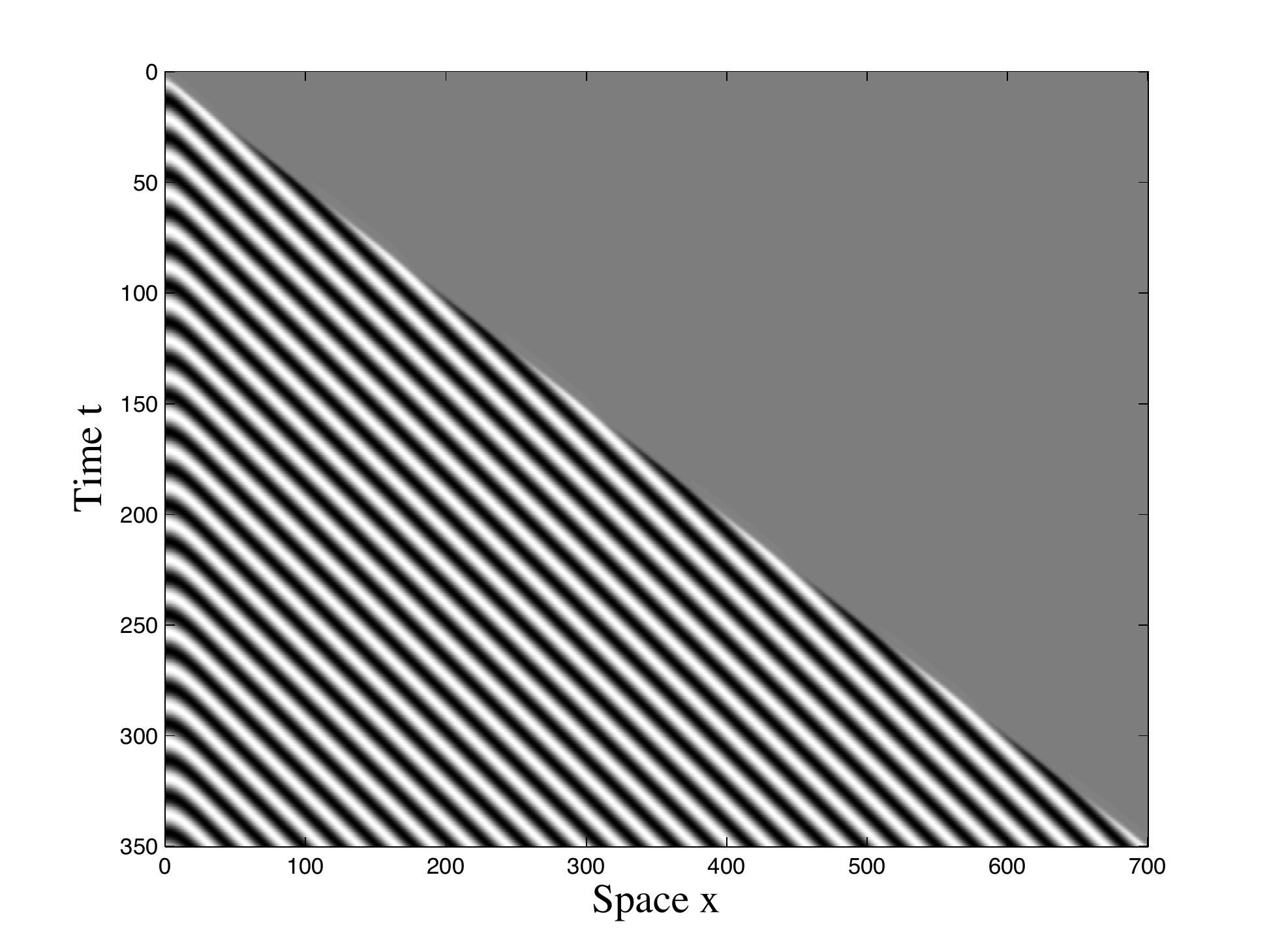}
        \end{subfigure}%
        \begin{subfigure}[h!]{0.45\textwidth}
                \centering
                \includegraphics[width=\textwidth]{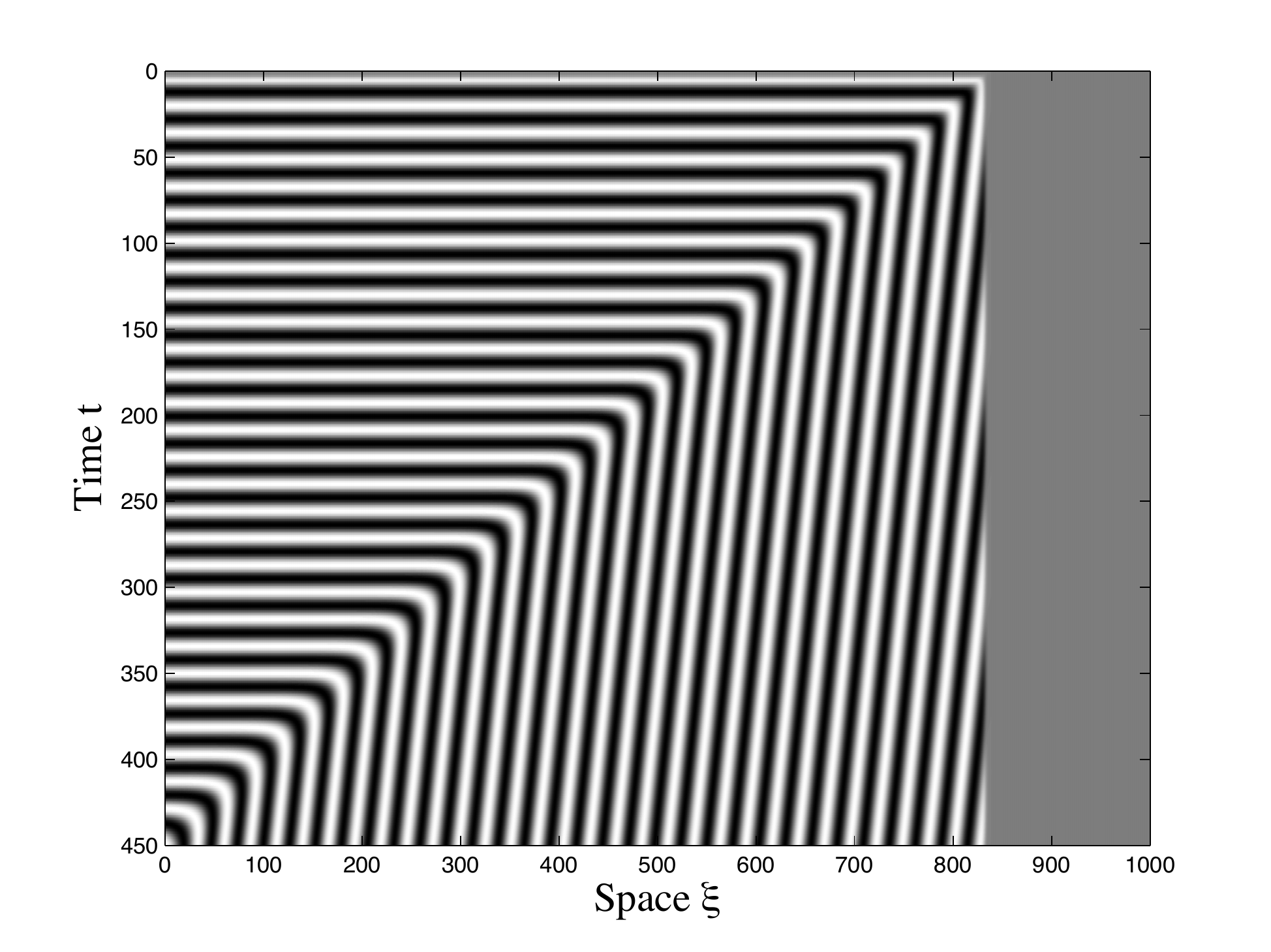}
        \end{subfigure}
        \vspace{-0.2in}
        \caption{Simulations of trigger fronts in both stationary (left) and moving (right) frames; $\,\alpha = -0.1, \gamma = -0.4, c = 1.8 <c_\rlin = 2\sqrt{1+\alpha^2} \approx 2.01.$ }\label{fig:spacetime}
\end{figure}

\begin{figure}[h!]
        \centering
        \begin{subfigure}[h!]{0.45\textwidth}
                \centering
                \includegraphics[width=\textwidth]{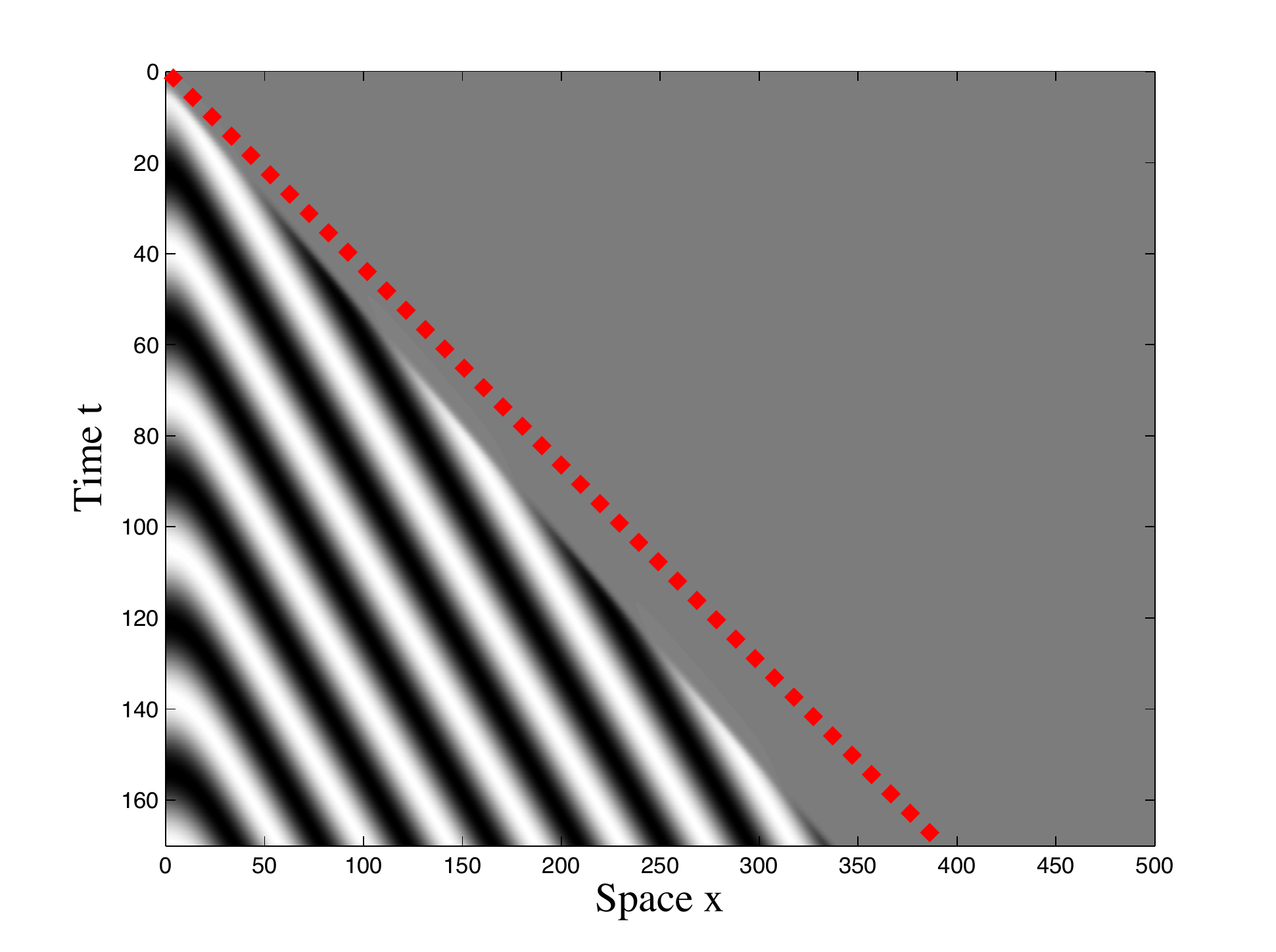}

        \end{subfigure}%
             \begin{subfigure}[h!]{0.45\textwidth}
                \centering
                \includegraphics[width=\textwidth]{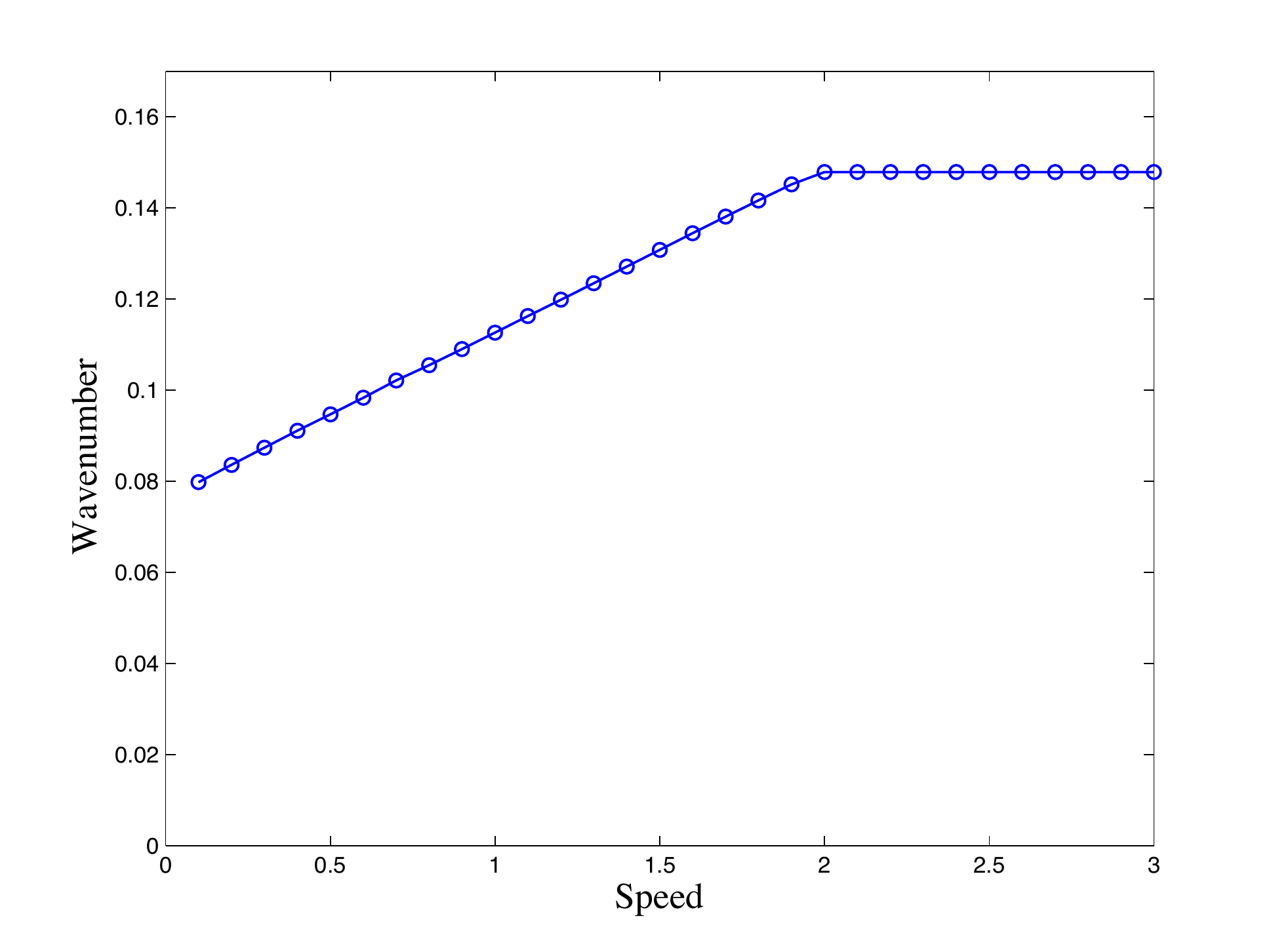}               \end{subfigure} 
              \vspace{-0.2in}
        \caption{The figure on the right gives the selected wavenumber for a range of trigger speeds;  $\alpha = -0.1,\gamma = -0.2$.  The figure on the left shows the space-time plot in the case $c = 2.3>c_\rlin $, when the constant wavenumber $k_\rlin$ is selected. The dotted red line denotes the path of the trigger $x = c t.$}\label{fig:2}
        \end{figure}

Phenomenologically, one observes roughly two different regimes.  For large speeds $c$, one observes patterns in the wake of a front that propagates with a speed $c_\mathrm{lin}<c$. The patterns created in the wake of that front are roughly independent of the speed $c$. When the trigger speed is decreased below $c_\mathrm{lin}$, patterns nucleate roughly at the location of transition from stability to instability. The wavenumber of patterns in the wake now depends smoothly on the speed $c$. At the transition, $c\sim c_\mathrm{lin}$, the wavenumber in the wake depends continuously on $c$, constant for $c>c_\mathrm{lin}$, and, to 
leading order, linear for $c\lesssim c_\mathrm{lin}$. This scenario is illustrated in space-time plots of direct numerical simulations in Figure \ref{fig:spacetime} and Figure \ref{fig:2}. There, the trigger moves into the medium from the left, sending out a periodic wave-train in it's wake.  The right-hand figure shows the same system in a moving frame of speed $c$.  The trigger is placed near the right boundary of the domain and sends out waves to the left. Since $A\equiv0$ is unstable for $\xi<0$, without any local perturbations, the system oscillates homogeneously.  As time evolves, a periodic wave train emitted from the trigger invades these homogeneous oscillations along a sink \cite{Scheel03,vanSaarloos92}.

\vspace{-0.2in}

\paragraph{Free fronts.}
Our results are based on perturbing from the \emph{freely invading front}, which is the front that one would observe when $\chi\equiv \chi_-$. In that case, the medium is translation invariant and unstable, and the instability spreads in space with a typical characteristic speed, $c_\mathrm{lin}$. This invasion process has received much attention in the literature and we refer to the review article by van Saarloos \cite{vanSaarloos03} for an exposition of basic theory and applications. In the termiology used there, free fronts in the complex Ginzburg-Landau equation are \emph{pulled fronts}, with speed $c_\mathrm{lin}$ determined by the linearization at the trivial state $A\equiv 0$: compactly supported initial conditions initially grow exponentially, \emph{pointwise}, in coordinate frames moving with speed $|c|<c_\mathrm{lin}$, and decay exponentially, pointwise, for speeds $|c|>c_\mathrm{lin}$. At the critical speed $c_\mathrm{lin}$, one observes algebraic decay, superimposed on oscillations with a 
frequency $\omega_\mathrm{lin}$. Both speed $c_\mathrm{lin}$ and frequency $\omega_\mathrm{lin}$ can be computed from the linear dispersion relation, obtained from inserting the ansatz $A\sim \re^{\lambda t + \nu x}$ into \eqref{eq:cgl}. In our case this relation takes the form $d(\lambda,\nu)=(1+\ri\alpha)\nu^2 +1-\lambda$, by solving\footnote{In general, one also needs to verify a pinching condition, which however in our case is automatically satisfied.}
\[
d(\ri\omega-c\nu,\nu)=0,\qquad \frac{\rmd}{\rmd\nu}d(\ri\omega-c\nu,\nu)=0
\]
where $\fr{d}{d\nu}$ is evaluated as a total derivative.
One finds 
\begin{equation}
c_\mathrm{lin}=2\sqrt{1+\alpha^2},\qquad \omega_\mathrm{lin}=\alpha, \qquad \nu_\mathrm{lin}=-\frac{1-\rmi\alpha}{\sqrt{1+\alpha^2}}.\label{e:linp}
\end{equation}
One derives a \emph{linear prediction for the nonlinear pattern} in the wake of the front from this linear information as follows. Nonlinear spatio-temporally periodic  patterns are, in the simplest form, solutions of the form 
\[
A(t,x)=\re^{\ri\Omega_\mathrm{nl}(k)t}A^\mathrm{p}(x;k),\qquad A^\mathrm{p}(x;k)=\sqrt{1-k^2}\re^{-\ri kx},
\quad \mbox{where }\,\, \Omega_\mathrm{nl}(k)=-\alpha k^2-\gamma(1-k^2).
\]
In the comoving frame of speed $c$, the frequency of these patterns is determined by the nonlinear dispersion relation in the comoving frame, 
\begin{equation}\label{eq:nldisp}
\Omega(k;c)=-\alpha k^2-\gamma(1-k^2)-ck.
\end{equation}
One can solve $\Omega(k,c_\mathrm{lin})=\omega_\mathrm{lin}$ for $k=:k_\mathrm{lin}$, using that $|k|<1$, and thereby derive a linear prediction for a nonlinear selected wavenumber,
\beq
k_\mathrm{lin} = \left\{\begin{array}{ll}\displaystyle{-\fr{\sqrt{1+\gamma^2}-\sqrt{1+\alpha^2}}{\gamma-\alpha}, \quad  }&\text{for} \ \gamma\neq\alpha\\
\displaystyle{-\fr{\alpha }{ \sqrt{1+\alpha^2}}, \quad }&\text{for}\ \gamma = \alpha.
\end{array}\right.
\label{eq:nlwn}
\eeq
We note that the group velocity of wave trains in the wake of the invasion front points away from the front interface,
\begin{equation}\label{e:lind}
c_\mathrm{g}:=\frac{\rmd\Omega(k;c_\mathrm{lin})}{\rmd k}\Big|_{k=k_\mathrm{lin}}=-2\sqrt{1+\gamma^2}<0.
\end{equation}
In other words, the invasion front acts as a wave source in its comoving frame; see \cite{Scheel03}. 

In order to state our main assumption, we first give a precise characterization of invasion fronts in the form needed later. 

\begin{Definition}[Generic free fronts]\label{d:1}
A free front is a solution to CGL, (\ref{e:cglt}) with $\chi \equiv 1$, of the form $A(t,x)=\rme^{\rmi\omega_\mathrm{lin} t} A^\mathrm{f}(x-c_\mathrm{lin}t)$, that satisfies
\[
A^\mathrm{f}(\xi)\to 0 \text{ for } \xi\to\infty,\qquad 
|A^\mathrm{f}(\xi)-A^\mathrm{p}(\xi;k_\mathrm{lin})|\to 0 \text{ for } \xi\to-\infty.
\]
We say a free front is generic if there exist $\mathcal{A}_\infty,\mathcal{B}_\infty \in \mb{C}$ with $\mathcal{A}_\infty\neq 0$ such that
\[
A^\mathrm{f}(\xi)=\left(\mathcal{A}_\infty \xi\rme^{\nu_\mathrm{lin}\xi}+\mathcal{B}_\infty\rme^{\nu_\mathrm{lin}\xi}\right)\Big(1+\rmo_\xi(1)\Big), \text{ for } \xi\to \infty.
\]
\end{Definition}

\paragraph{Trigger fronts and main result.}
Our goal is to find solutions to (\ref{e:cglt}) that describe the triggered invasion process when the speed of the trigger is less than, but close to, the speed of the invasion front. In a completely analogous fashion to free fronts, we define trigger fronts as spatial connections between $A=0$ at $\xi=+\infty$ and a periodic wave train at $\xi=-\infty$.

\begin{Definition}[Trigger fronts]
A trigger front with frequency $\omega_\mathrm{tf}$ is a solution to CGL, (\ref{e:cglt}), with trigger speed $c$, of the form $A(t,x)=\rme^{\rmi\omega_\mathrm{tf} t} A^\mathrm{tf}(x-ct)$, that satisfies
\[
A^\mathrm{tf}(\xi)\to 0 \text{ for } \xi\to\infty,\qquad 
|A^\mathrm{tf}(\xi)-A^\mathrm{p}(\xi;k_\mathrm{tf})|\to 0 \text{ for } \xi\to-\infty,
\]
where $k_\mathrm{tf}$ is such that $\omega_\mathrm{tf}=\Omega(k_\mathrm{tf};c)$ and the group velocity
\[
c_\mathrm{g}(k_\mathrm{tf}):=\frac{\rmd\Omega(k;c)}{\rmd k}\Big|_{k=k_\mathrm{tf}}<0.
\]
In other words, trigger fronts are time-periodic with frequency $\omega_\mathrm{tf}$ and emit wave trains with wavenumber $k_\mathrm{tf}$. Furthermore, for a fixed $\delta>0$ sufficiently small, we define the front interface distance as $$\xi_*:=\inf\left\{\xi;\sup_{\xi'>\xi} |A(\xi')|<\delta\right\}$$.
\end{Definition}

We are now ready to state our main result.

\begin{Theorem}\label{t:1} Fix $\alpha,\gamma\in \mb{R}$ and assume that there exists a generic free front. Then there exist trigger fronts for $c<c_\mathrm{lin}$, $|c-c_\mathrm{lin}|$ sufficiently small. The frequency of the trigger front possesses the expansion
\beq\label{eq:omsel}
\omega_\mathrm{tf}( c) =  \omega_\rab(c) + \fr{2 |\Delta Z_\ri|}{\pi (1+\alpha^2)^{3/4}}  ( c_\rlin -  c)^{3/2} + \rmO( ( c_\rlin -  c)^{2}).
\eeq
Here, 
\[
\omega_\mathrm{abs}(c)=-\alpha+\frac{\alpha c^2}{2(1+\alpha^2)},
\]
and $\Delta Z_\ri$ is defined in (\ref{def:DI}), below.
Furthermore, for $\alpha \neq \gamma$ the selected wavenumber has the expansion
\[
k_\mathrm{tf} = k_\rlin - g_1(\alpha,\gamma) (c_\rlin - c)  -   \fr{2 |\Delta Z_\ri|}{\pi (1+\gamma^2)^{1/2}(1+\alpha^2)^{3/4}}  ( c_\rlin -  c)^{3/2} + \mc{O}(( c_\rlin -  c)^{2})
\]
where $g_1(\alpha,\gamma) =\fr{1}{2(\gamma - \alpha)}\lp(\fr{1+2\alpha\gamma-\alpha^2 }{\sqrt{(1+\alpha^2)(1+\gamma^2)}}   +1\rp).$ 
The distance between the trigger and the front interface is given by
\[
{\xi}_* = \pi(1+\alpha^2)^{3/4} ( c_\rlin -  c)^{-1/2} + (1+\alpha^2)^{1/2}\, \Delta Z_\mathrm{r} + \rmO(( c_\rlin -  c)^{1/2}),
\]
where $\Delta Z_\mathrm{r}$ is defined in (\ref{def:DI}) as well.
\end{Theorem}

In the following, we elaborate on assumptions and conclusions in this result.

\begin{figure}[h!]
\centering

\begin{subfigure}[h!]{0.48\textwidth}
  \def\svgwidth{200pt}
  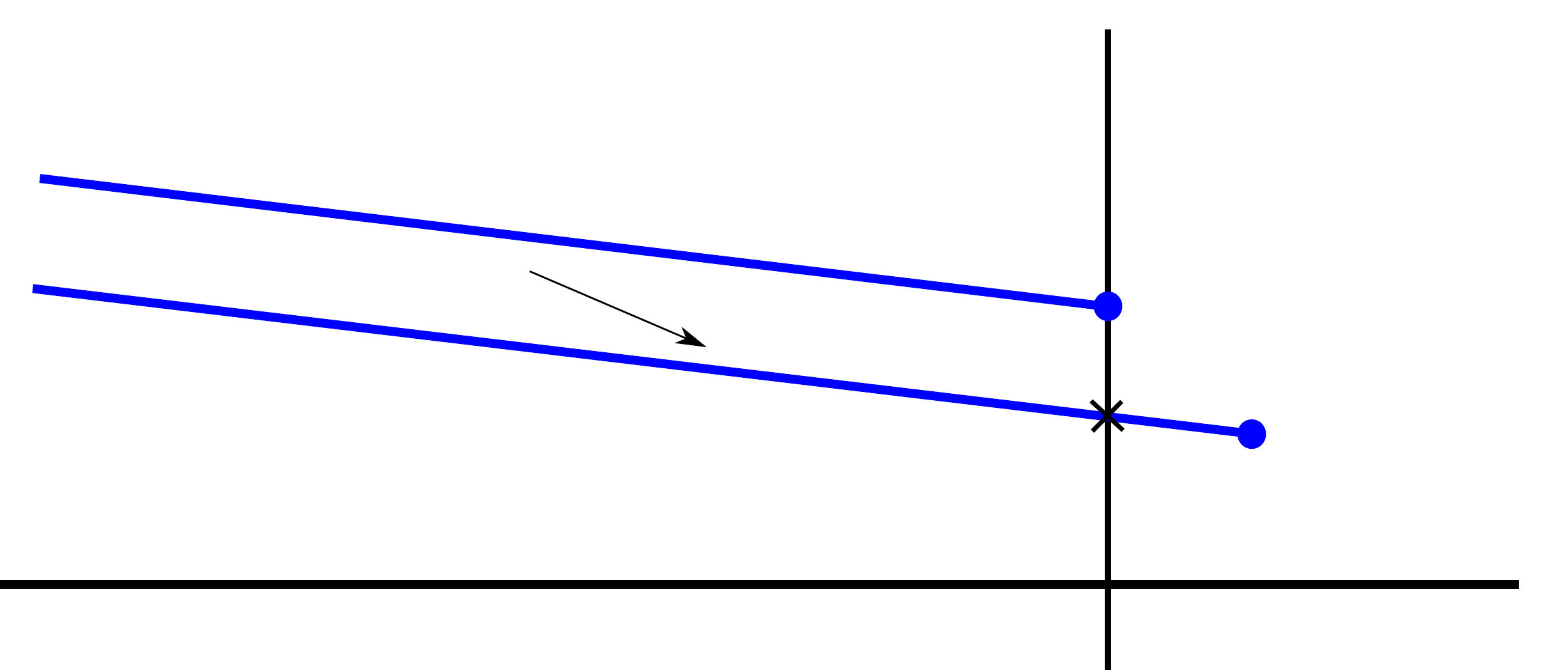
  \end{subfigure}
  \vspace{-0.1in}
  \begin{subfigure}[h!]{0.48\textwidth}
  \centering
    \def\svgwidth{200pt}
  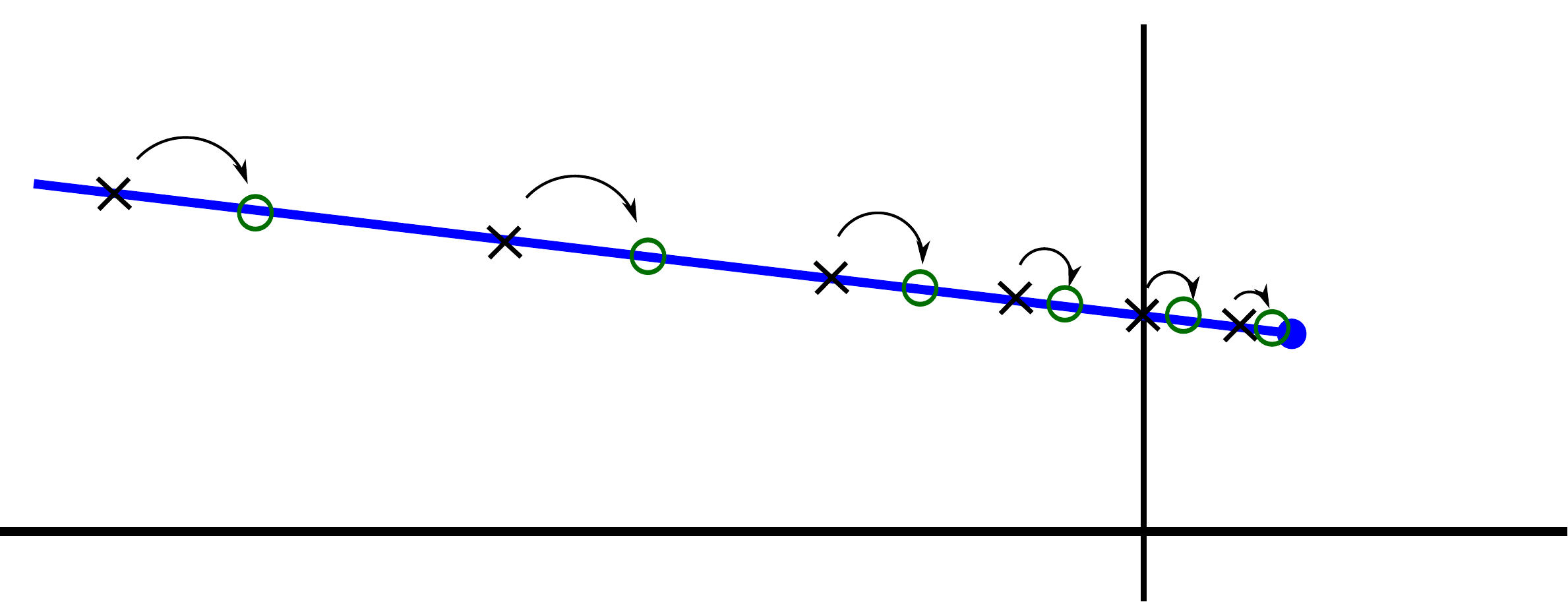
  \end{subfigure}
  
\caption{Left: Plot of the absolute spectrum $\Sigma_\rab$ in the complex plane. The blue dots denote double roots. The curve moves into the right half plane for $c<c_\rlin$, intersecting $\ri\R$ at $\ri\omega_\rab$. Right: As the domain size $L$ increases to $+\infty$ the spectrum of the linear operator accumulates on $\Sigma_\rab$. }
\label{fig:absspec}
\end{figure}

\paragraph{Interpreting the expansion.}
The frequency $\omega_\mathrm{abs}$ is determined by the intersection of the absolute spectrum \cite{Scheel00} of the linearization at the origin in a frame moving with speed $c$ and the imaginary axis. Roughly speaking, the absolute spectrum denotes curves in the complex plane that emanate from double roots of the dispersion relation such that finitely truncated boundary-value-problems possess dense clusters of eigenvalues at these curves as the size of the domain goes to infinity \cite{Scheel00}. The absolute spectrum is a natural first-order prediction for frequencies of trigger fronts as we shall explain in Section \ref{sec:dynbp}. Intersection points with the imaginary axis arise when the edge of the absolute spectrum crosses the imaginary axis. This happens precisely when $c$ is decreased below $c_\mathrm{lin}$, so that quite generally there will be a unique intersection point, smoothly depending on $c$;
see Figure \ref{fig:absspec}.

The term $\Delta Z_\ri$ that describes corrections to the leading-order prediction can be interpreted as a distance in projective space between tangent spaces of stable manifold and unstable manifold. Roughly speaking, decay at $\xi=+\infty$ creates an effective boundary condition of Robin type at $\xi=0$, while the leading edge of the invasion front approximately satisfies a different Robin boundary condition. The distance between those two boundary conditions, measured in an appropriate coordinate system, is denoted by $\Delta Z_\ri$. Using a simple shooting algorithm, one can evaluate $\Delta Z_\ri$ quite accurately; see our numerical studies in Section \ref{sec:num}.

\paragraph{The assumption on existence, non-uniqueness, and  more fronts.}
We will see in our proof that the expansion of invasion fronts at $\xi\to\infty$ holds for any invasion front with $|\mathcal{A}_\infty|+|\mathcal{B}_\infty|>0$, so that genericity refers to the (open) condition $\mathcal{A}_\infty\neq 0$, only. Existence can be cast as the problem of existence of a heteroclinic connection between an unstable equilibrium and a sink. Genericity refers to unstable and stable manifolds being in general position. 

We are not aware of results guaranteeing the existence of generic invasion fronts, or any evidence pointing towards non-existence. We show in Proposition \ref{prop:genhet} that generic invasion fronts exist for $|\alpha-\gamma|$ sufficiently small. 

In a different interpretation, vanishing of $\mathcal{A}_\infty$ characterizes fronts at the boundary between pushed and pulled invasion in parameter space. In the case of cubic CGL, considered here, such a transition has not been observed. On the other hand, the transition is usually accompanied by an increase in the (nonlinear) invasion speed, so that one expects trigger fronts to exist for speeds larger than $c_\mathrm{lin}$ beyond this transition; see the discussion in Section \ref{sec:disc}.

The fronts we find are \emph{not unique}. In fact, our proof gives the existence of a countable family of invasion fronts, indexed by $j=1,2,\ldots$, with frequencies 
\beq\label{eq:omselj}
\omega_\mathrm{tf}( c;j) =  \omega_\rab(c) + \fr{2 |\Delta Z_\ri|}{\pi j (1+\alpha^2)^{3/4}}  ( c_\rlin -  c)^{3/2} + \rmO( ( c_\rlin -  c)^{5/2}).
\eeq
Roughly speaking, the distance between front interface (measured by, say, the location of $|A|=\delta>0$, fixed) and the trigger location $\xi=0$ increases linearly with the index $j$. We expect these fronts to be unstable with Morse index increasing linearly in $j$.

\paragraph{Stability and secondary instabilities.}

\begin{figure}[h!] 
\centering
   \begin{subfigure}[h!]{0.45\textwidth}
     \includegraphics[width=\textwidth]{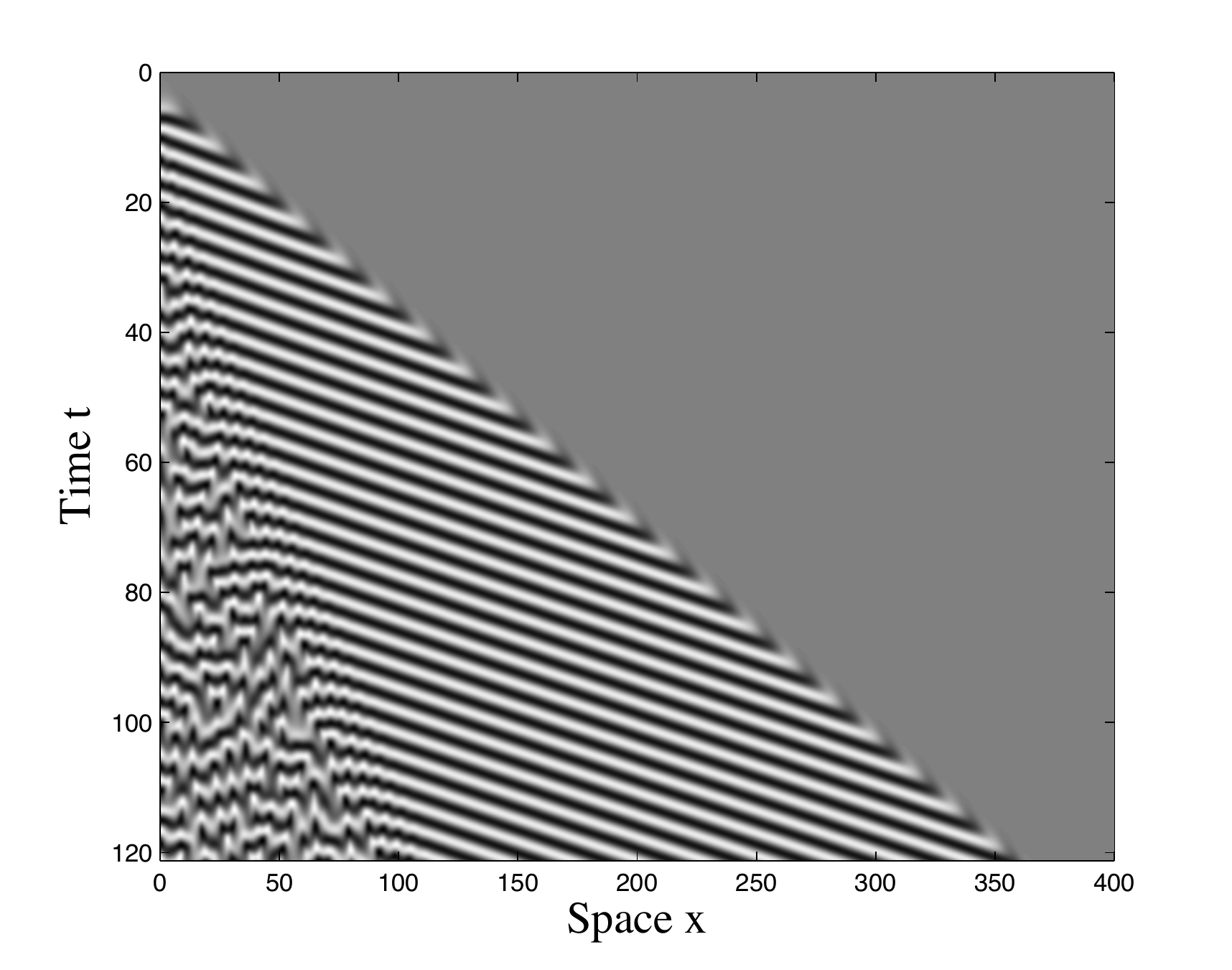}
   \end{subfigure}%
   \begin{subfigure}[h!]{0.45\textwidth}
     \includegraphics[width=\textwidth]{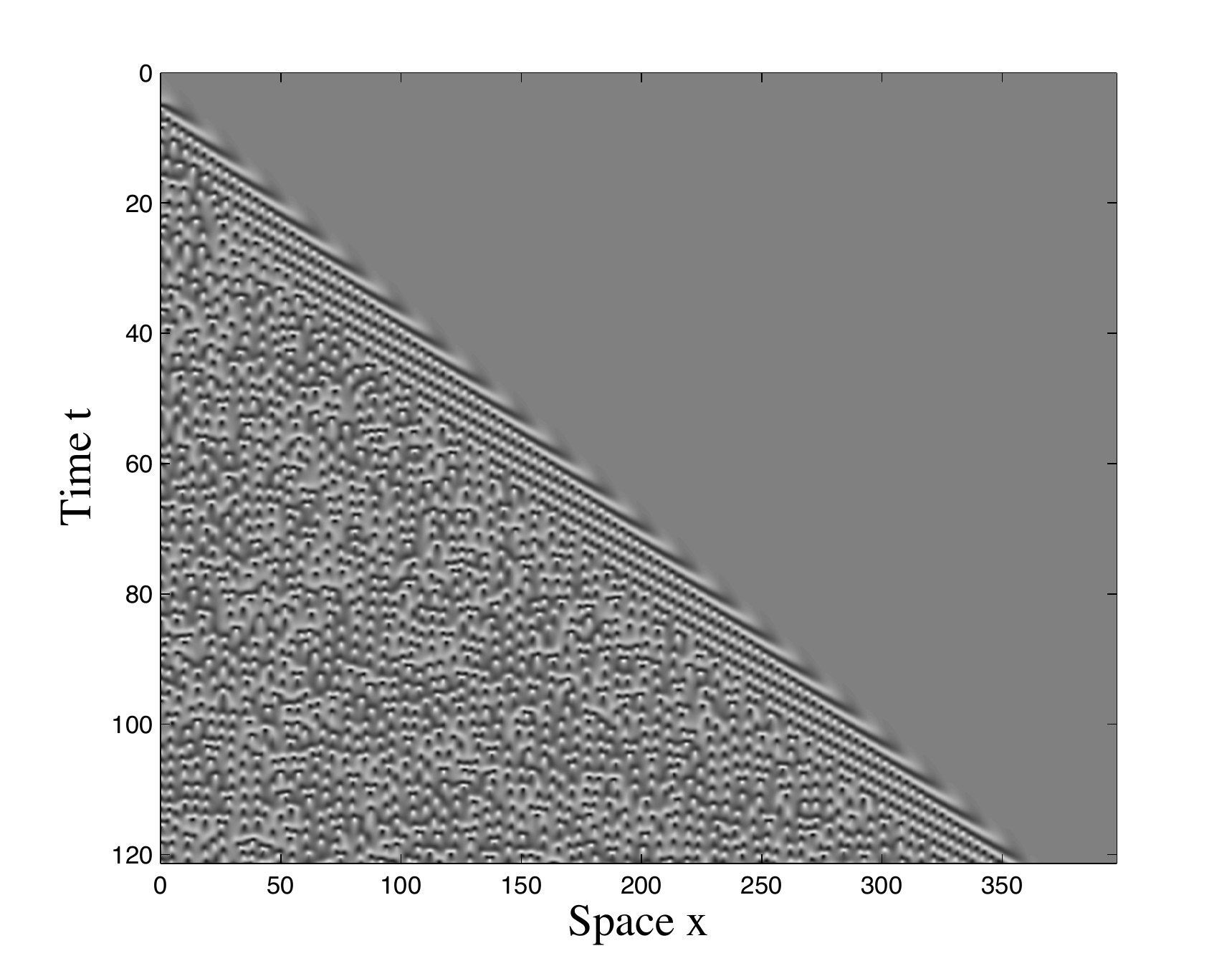}
   \end{subfigure}
\caption{Invasion of spatio-temporal chaos after formation of periodic pattern in wake of a trigger; parameters in the Benjamin Feir instability range: $\alpha = -1.2,c = 0.95\cdot  c_\rlin = 1.90\cdot \sqrt{1+\alpha^2}$ and $\gamma = 2$ (left), and $\gamma = 8$ (right).} \label{fig:ch1}
\end{figure}

We did not attempt to prove stability but we suspect that the fronts we find are stable in a suitable sense. In fact, stability of the free front is known for $\alpha=\gamma=0$ \cite{eckwayne} and one can show that fronts with $\alpha\sim\gamma$ are at least linearly stable. It would be interesting to conclude stability of the trigger fronts that we find from spectral stability of the free front. On the other hand, free fronts can be unstable. In particular, the wave train selected by the invasion (or trigger) front is in fact unstable for many values of $\alpha$ and $\gamma$. Most dramatically, for $\alpha\gamma<-1$, \emph{all} wave trains are unstable. The instability may propagate slower than the primary invasion front, leading to  a wedge in space-time plots where the (unstable) selected wave train can be observed.  For large $\alpha,\gamma$, the secondary instability invades at a fixed distance behind the primary front; see \cite[\S 2.11.5]{vanSaarloos03}, \cite{radesherr} and Figure \ref{fig:ch1} for an illustration.

\paragraph{Spatial dynamics.}

We prove our main result by rewriting the existence problem for a trigger front as a non-autonomous ODE for $A^\mathrm{tf}(\xi)$,
\begin{equation}
A_{\xi\xi} = - \fr{1}{(1 + \ri {\alpha})} \lp[ (\chi(\xi) -  i {\omega}) A + {c} A_\xi  - (1 + {\gamma} \ri) A |A|^2\rp]\label{eq:us-mf1}.
\end{equation}
Here, $c$ is the (externally prescribed) bifurcation parameter and $\omega$ is a matching parameter that will be used to achieve appropriate intersections of stable and unstable manifolds. In fact, one can most easily understand this ODE by separately considering dynamics with $\chi\equiv \chi_+$ and $\chi=\chi_-$, thus obtaining two separate 4-dimensional phase portraits. For $\chi=\chi_+$, we will find the equilibrium $A=A_\xi=0$ with a two-dimensional stable manifold. For $\chi=\chi_-$, we find that  $A=A_\xi=0$ is stable with an open basin of attraction. On the other hand, we also have a periodic orbit $A^\mathrm{p}$, depending on $c$ and $\omega$, with a two-dimensional unstable manifold. Existence of a free front implies that there exists a heteroclinic orbit between $A^\mathrm{p}$ and the origin. The strategy of the proof is to find intersections between the unstable manifold of $A^\mathrm{p}$ in the $\chi_-$-phase-space and the stable manifold of the origin in the $\chi_+$-phase space. We will find 
these intersections bifurcating from the origin upon decreasing $c$ below $c_\mathrm{lin}$ and adjusting $\omega$ appropriately. A key technical ingredient to our result is a geometric blowup of the origin which both factors the $S^1$-action on $(A,A_\xi)\in \C^2$ induced by the gauge symmetry, and desingularizes the vector field by effectively separating eigenspaces near an algebraically double complex eigenvalue. Intersections can then be found almost explicitly on a singular sphere, and lifted using Fenichel's invariant foliation methods \cite{Fenichel71,Fenichel74,Fenichel77}.

\paragraph{Outline.}
The remainder of this paper is organized as follows.
 Section \ref{sec:dynbp} gives heuristics and a  conceptual outline of the proof of Theorem \ref{t:1}, describing the role of the absolute spectrum in finding the desired heteroclinic intersection. 
Section \ref{sec:pre} contains the proof of our main theorem. 
Section \ref{sec:num} gives comparisons between our expansions, direct simulations, and direct computations of heteroclinic orbits.
Finally, Section \ref{sec:disc} gives a discussion of our results and directions for future research.

  \begin{Acknowledgment}
A. Scheel was partially supported by the National Science
Foundation through grant NSF-DMS-0806614. This material is based upon work supported by the National Science Foundation Graduate Research Fellowship under grant NSF-GFRP-00006595. Any opinion, findings, and conclusions or recommendations expressed in this material are those of the authors(s) and do not necessarily reflect the views of the National Science Foundation.
\end{Acknowledgment}

\section{Heuristics --- formal asymptotics and the role of the absolute spectrum.}\label{sec:dynbp}
We present several concepts that help understand the result. We first explain the role of absolute spectra and then present a formal \emph{exponential matching} argument that mimicks the main strategy of proof. Last, we give a brief outline of \emph{geometric desingularization}, the center piece of our approach and connect it with these two other points of view.

\paragraph{Absolute spectra.} Taking the perspective of a numerical simulation in a comoving frame, we could attempt to understand trigger fronts as steady states (up to the gauge symmetry) bifurcating from the trivial state $A\equiv 0$ when the latter looses stability. One would then linearize the system at $A\equiv 0$ in a bounded domain, in a comoving frame, equipped with separated boundary conditions. In order to realistically capture the phenomena, we would assume that the size of the domain $L$ is large. One can even simplify further and substitute ``effective'' boundary conditions at $\xi=0$ and restrict to $\xi\in[-L,0]$. The linearized operator $(1+\rmi\alpha)\partial_{\xi\xi}+c\partial_\xi+1-\rmi\omega$ then possesses constant coefficients but is not self-adjoint and spectra in large domains are not approximated by spectra in unbounded domains. In fact, eigenvalues of the linearized operator in the bounded interval accumulate, as $L\to\infty$, at curves referred to as the \emph{absolute spectrum}, $\Sigma_\mathrm{abs}$ \cite{Scheel00}; see Fig. \ref{fig:absspec}. Those curves can be computed from the dispersion relation $d(\lambda,\nu)=0$ as follows. One varies $\lambda$ as a parameter, solves for all roots $\nu=\nu_j(\lambda)$, and orders those roots by real part $\Re\,\nu_1\leq \Re\,\nu_2\leq\ldots \leq\Re\,\nu_N$ for each fixed $\lambda$. For $\Re\,\lambda\to\infty$, one always finds $\Re\,\nu_j<0<\Re\,\nu_{j'}$, for $j\leq M<j'$, with some fixed $M$. The absolute spectrum then is the set of $\lambda$ so that $\Re\,\nu_M(\lambda)=\Re\,\nu_{M+1}(\lambda)$.

In our case, $d_c(\lambda,\nu)=(1+\rmi\alpha)\nu^2+c\nu+1-\lambda$, with roots
$$
\nu_{\pm} = -\fr{ c}{2(1+\ri \alpha)} \pm\sqrt{\fr{ c^2}{4(1+\ri\alpha)^2} - \fr{1 - \lambda}{1+\ri\alpha}}.
$$
At the absolute spectrum, we must have $\rre\,\nu_+ = \rre\,\nu_-$, so  
\begin{equation}\label{e:discr}
\fr{ c^2}{4(1+\ri\alpha)^2} - \fr{1 - \lambda}{1+\ri\alpha} <0.
\end{equation}
Bifurcations in a bounded domain occur when eigenvalues cross the imaginary axis, $\lambda=\rmi\omega$, leading to periodic orbits with frequency close to $\omega$. We are therefore interested in the intersection of the absolute spectrum with the imaginary axis, which, starting with (\ref{e:discr}), is readily found at
\beq 
\omega_\rab(c) = -\alpha+\fr{ \alpha  c^2}{2(1+\alpha^2)}, 
\eeq
for $c<c_\mathrm{lin}$. 
We caution the reader that $\omega_\mathrm{abs}$ is not constant, due to two effects: the curve of absolute spectrum is not horizontal and the imaginary part of the leading edge will depend on $c$.
As the domain size increases, eigenvalues will move along the curve of absolute spectrum towards the edge \cite[Lemma 5.5]{Scheel00}, which is located in $\{\Re\lambda>0\}$, thus leading to a sequence of Hopf bifurcations. The periodic solutions with frequency $\omega$ emerging from these Hopf bifurcations converge to the trigger fronts described in our main result; see also the subsequent discussion, pointing to a countable family of trigger fronts. 

While this view point gives a rather simple and universal prediction, it is generally not clear how the bifurcating solutions evolve as $L\to\infty$. One can easily envision pushed fronts leading to a faster propagation mechanism, thus inducing turning points in the branch of periodic solutions. We refer to \cite{couairon,Tobias98,basins} for more detailed analyses in this direction. 

Inspecting the proofs in \cite{Scheel00}, eigenvalues are generated by oscillatory dynamics in the Grassmannian: one evolves the boundary condition at $\xi=-L$ as an element of the Grassmannian and seeks intersections with the boundary condition at $\xi=L$. Since eigenvalues in the Grassmannian are differences of eigenvalues in the linear system, oscillatory dynamics occur when $\Re\nu_-=\Re\nu_+$, in our simple case. In the following, we explain how these oscillatory dynamics can be found in a more local analysis.

\paragraph{Exponential matching.} Starting with the assumption that the amplitude of the trigger front would be small at $\xi=0$ when $c\lesssim c_\mathrm{lin}$, we try to glue the free front, shifted appropriately, at $\xi=0$ with an exponentially decaying tail in $\xi>0$. 

To leading order, the free front decays like $\mathcal{A}_\infty\zeta\rme^{\nu(\xi+\zeta)}+\mathcal{B}_\infty\rme^{\nu(\xi+\zeta)}$, with shifted variable $\zeta$. We next vary $\omega$ and $c$ and neglect the dependence of the coefficients $\mathcal{A}_\infty$ and $\mathcal{B}_\infty$ on $c,\omega$. In order to track the dependence of the exponential rate $\nu$ on $\omega,c$, we inspect the linear dispersion relation $(1+\ri\alpha)\nu^2 + c\nu +1-\rmi\omega=0$. Near $\omega=\omega_\mathrm{lin}$ and $c=c_\mathrm{lin}$, we have two roots $\nu_1\sim\nu_2$, so that one expects decay in the front profile $\mathcal{A}\rme^{\nu_2(\xi+\zeta)}$, when $\Re\,\nu_1<\Re\,\nu_2<0$. 

On the other side, small solutions in $\xi>0$ decay like $\mathcal{A}_+\rme^{\nu_+\xi}$, where $\nu_+$ solves $(1+\ri\alpha)\nu^2 + c\nu -1-\rmi\omega=0$. 

Trying to match $A$ and $A_\xi$ at $\xi=0$, we find that due to leading order scaling invariance (keeping only linear terms) and gauge invariance, it is enough to match $A_\xi/A$. This however gives $\nu_+=\nu_1$ in our leading order expansion, an equation that cannot be solved varying $\omega$ and $c$ locally. In fact, our expansion for the invasion front $\mathcal{A}\rme^{\nu_1(\xi+\zeta)}$ is not valid uniformly when varying $\omega$ and $c$. In order to get smooth and uniform expansions, one needs to use a smooth expansion in terms of exponentials near the double root, using either smooth normal forms near a Jordan block \cite{arnold} or, more along the lines of our strategy, normal forms for the induced flow on $A_\xi/A$. Without going into details, one readily finds that the expansion in $\rme^{\nu_2\xi}$ alone will fail near $\Re\,\nu_1=\Re\,\nu_2$, where one would keep both terms $\mathcal{A}\rme^{\nu_1\zeta}+\mathcal{B}\rme^{\nu_2\zeta}$. Varying $\zeta$ and the difference $\nu_1-\nu_2$ as a function of 
$\omega$, one can then solve the matching problem with $\xi>0$. 

Note that here the absolute spectrum appears in a natural way as parameter values $\omega$ where $\Re\nu_1=\Re\nu_2$. 

\paragraph{Geometric desingularization.}

\begin{figure}[htb]
  \centering
  \def\svgwidth{300pt}
  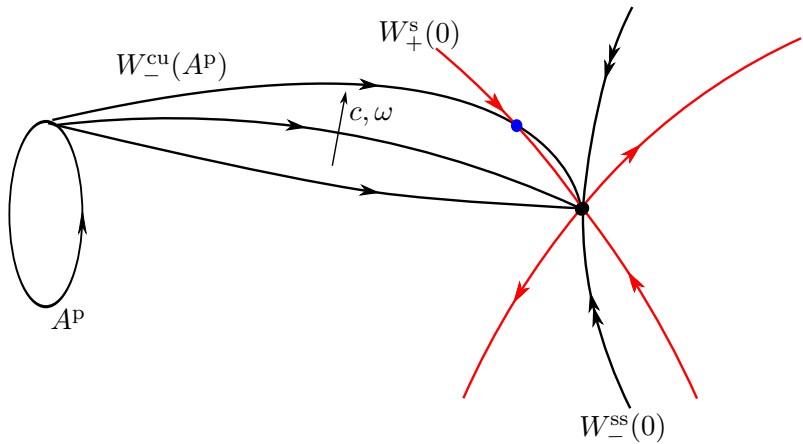
  \caption{Schematic diagram of intersection between $W_-^\rcu(A^\rpp)$ and $W_+^\rs(0)$.  The red lines denote the invariant manifolds of the origin in the $\xi>0$ dynamics.  These act as boundary conditions for a shooting from $W_-^\rcu(A^\rpp)$.  Unfolding $c$ and $\omega$ from $c_\rlin$ and $\omega_\rlin$, we obtain a non-trivial intersection at $(\xi, \omega) = (\xi_*,\omega_\rtf(c))$  (blue dot)
  .}\label{dia:mflds}
\end{figure}

The global aspect of the matching procedure is illustrated in Figure \ref{dia:mflds}. The dynamics with $\chi\equiv \chi_-$ show a periodic orbit $A^\mathrm{p}$ with a two-dimensional unstable manifold that converges to the stable equilibrium $A\equiv 0$. In this phase portrait, the linearization at the origin is a complex Jordan block. Overlayed is the phase portrait for $\chi\equiv \chi_+$, where $A\equiv 0$ possesses a two-dimensional stable manifold. Factoring the $S^1$-symmetry, both manifolds are 1-dimensional in a 3-dimensional ambient space. We are looking for intersections close to the origin. Our main control parameter is $\omega$, which, together with the bifurcation parameter $c$ unfolds the complex Jordan block and leads to possible flips in the position of the unstable manifold of $A^\mathrm{p}$ near $A\equiv 0$.

From the above discussion, it is apparent that a description of a neighborhood of $A=A_\xi=0$ is crucial to the understanding of the matching procedure. There are several dynamical systems techniques for such a description. First, one can try to use smooth linearization for the ODE (\ref{eq:us-mf1}). Since $\nu=\nu_\mathrm{lin}$ has negative real part, the double eigenvalue is in fact non-resonant, and smooth linearization results are available; see for instance \cite[\S 6.6]{katok}. In our context, we would however also require smooth parameter-dependence and $S^1$-equivariance. While such results may well be true, our approach appears more robust and elementary. Roughly speaking, we introduce polar coordinates for $(A,A_\xi)\in\C^2\sim \R_+\times S^3$, and factor the gauge symmetry, $\C^2/S^1\sim \R_+\times (S^3/S^1)\sim \R^+\times S^2$, where the last equivalence collapses the Hopf fibration. Identifying $S^2$ with the Riemann sphere $\overline{\mb{C}}$, we end up with coordinates $|A|^2$ and $A_\xi/A$, that we 
identified above. In these coordinates, the sphere $|A|^2=0$ is normally hyperbolic and carries an explicit linear projective flow that allows us to compute expansions at leading order. Matching outside of the sphere, which corresponds to taking higher-order terms in the exponential matching procedure into account, can be readily achieved after straightening out smooth foliations. We refer to the next sections for details of this procedure. 

\section{Heteroclinic bifurcation analysis}\label{sec:pre}

In this section we prove our main result. We start in Section \ref{s:3.1} by several simple rescaling transformations and calculate dimensions of stable and unstable manifolds using information from the dispersion relation. Section \ref{s:3.3} introduces the coordinate system near the origin that we use to factor the gauge symmetry and desingularize the linear dynamics. Section \ref{s:3.4} examines the asymptotics of the free invasion front and shows that our assumptions are verified in the regime $\alpha\sim\gamma$. Section \ref{s:3.5} contains the matching analysis in the desingularized coordinates.

\subsection{Scalings and dimension counting}
\label{s:3.1}
The following change of variables will slightly simplify the remainder of our analysis, essentially eliminating $\alpha$ as a parameter. 
Suppose first that $1+\alpha\gamma>0$. Setting 
\beq\label{e:scal}
S= \fr{{c}{\alpha}}{2\sqrt{1+{\alpha}^2}},\quad m^2= 1 + \fr{( c  \alpha)^2}{4(1+\alpha^2)} - \alpha  \omega,\quad l^2= \fr{1+ \alpha  \gamma}{m^2},
\eeq
we scale and parameterize
\beq\label{e:scal2}
a = l \re^{-\ri S \xi} A,\quad 
\hat{c} = \fr{c}{m\sqrt{1 + \alpha^2}}, \quad  
\hat\omega =\fr{\omega  - \omega_\rab(\alpha, c)}{m^2}, \quad
\zeta = \fr{m}{\sqrt{1+ \alpha^2}}\xi  ,\quad 
\hat{\gamma} = \fr{\gamma -\alpha }{1 + \gamma\alpha},
\eeq
so that \eqref{eq:us-mf1} simplifies to 
\beq 
 a_{\zeta\zeta}= -(1 -  \ri\hat\omega) a - \hat{c}  a_\zeta+ ( 1+ \ri\hat\gamma) a |a|^2 + \fr{(1 - \chi)}{m^2}(1+\ri\alpha) a.\label{s-mfa}
\eeq
If we restrict to $\zeta<0$ we have $\chi = 1$ and obtain
\beq
a_{\zeta\zeta} = -(1 - \ri\hat\omega ) a - \hat{c}  a_\zeta+ ( 1+ \ri\hat\gamma) a |a|^2,\label{s-mfl}
\eeq
whereas for $\zeta\geq0$ we obtain 
\beq 
 a_{\zeta\zeta} = -(1 - \ri\hat\omega ) a - \hat{c}  a_\zeta+ ( 1+ \ri\hat\gamma ) a |a|^2 + \fr{2}{m^2}(1+\ri\alpha) a.\label{s-mf}
\eeq
The cases $1+\alpha \gamma< 0$ and $1+\alpha\gamma=0$ can be treated in a similar fashion. For $1+\alpha \gamma< 0$ one finds 
\beq
 a_{\zeta\zeta} = -(1 - \ri\hat\omega ) a - \hat{c}  a_\zeta+ ( -1+ \ri\hat\gamma ) a |a|^2 + \fr{(1 - \chi)}{m^2}(1+\ri\alpha) a\label{s-mf2}
,\eeq
and for $1+\alpha\gamma=0$ we find
\beq
 a_{\zeta\zeta} = -(1 - \ri\hat\omega ) a - \hat{c}  a_\zeta+  \ri\hat\gamma  a |a|^2 + \fr{(1 - \chi)}{m^2}(1+\ri\alpha) a\label{s-mf3}
.\eeq
Since those cases do not alter the subsequent discussion, we omit the straightforward details and focus on the case $1+\alpha\gamma>0$ in the sequel. 
\begin{Remark}
The fact that the scaled equation changes type when $1+\alpha\gamma$ changes sign reflects the Benjamin-Feir instability: for $1+\alpha\gamma<0$, all wave trains are unstable. One observes secondary invasion fronts, Figure \ref{fig:ch1}, which for moderate values of $|\alpha\gamma|$ propagate slower than the primary invasion front, leaving a long plateau of unstable wave trains. For large values of $|\alpha\gamma|$, the instability catches up with the primary front and one observes chaotic dynamics in the immediate wake of the primary front. In this case, our results, while correct, do not describe actually observed wavenumbers due to the strong, absolute instability in the wake; we refer to \cite[\S 2.11.5]{vanSaarloos03} and \cite{nozakibekki} for a discussion of this behavior for free invasion fronts.  
\end{Remark}
For the remainder of the proof, we will consider the first order form of (\ref{s-mfa}), writing $\zeta$ as our spatial variable, $A$ once again as the amplitude, and $B := A_\zeta$.
\begin{align}
A_\zeta &= B\notag\\
B_\zeta &=  -(1 -  \ri\hat\omega) A - \hat{c}  A_{\zeta}+ ( 1+ \ri\hat\gamma) A |A|^2 + \fr{(1 - \chi)}{m^2}(1+\ri\alpha) A.\label{s-mfa-n}
\end{align}


The following proposition shows that free fronts are unique up to spatial translation and gauge symmetry. 

\begin{Proposition}\label{lem:dim2}
Let $c = c_\rlin$ and $\omega = \omega_\rlin$.  Then the CGL traveling wave equation \eqref{eq:us-mf1} possesses a unique relative equilibrium $A^\mathrm{p}$ with wavenumber $k_\mathrm{lin}$. In other words, there exists a unique wave train with frequency $\omega_\rlin$ in a frame moving with speed $c_\rlin$. In \eqref{eq:us-mf1}, this relative equilibrium possesses a smooth two-dimensional center-unstable manifold $W_-^\rcu(A^\rpp)$.
\end{Proposition}
\begin{Proof}
The proof of the proposition is a direct calculation of eigenvalues of the linearization at $A^\mathrm{p}$. We refer to  \cite[\S 2.2.3]{vanSaarloos92}, where explicit conditions are given so that this linearization possesses precisely one unstable eigenvalue but caution the reader that the convention in this paper is to define frequency and wavenumber of wavetrains via $A = r \re^{-\ri \omega t + k x}$ (i.e. with the opposite sign as ours). Substituting $ c = c_\rlin$ and $\omega = \omega_\rlin$ into the formulas there, we obtain the desired result for all $\alpha,\gamma$.  
\end{Proof}

\subsection{Symmetry reduction and geometric blowup}
\label{s:3.3}

In order to carry out the matching procedure near the origin, we would like to quotient the $S^1$-action and exhibit the leading-order scaling symmetry from the linearized equation. It turns out that this can be achieved in a very simple fashion, introducing $|A|^2\geq 0$ and $A_\zeta/A\in\C$  as new variables. While very effective, this choice appears somewhat arbitrary and we will show how to obtain these coordinates in a systematic fashion.

Since the $S^1$-action is not free near the origin, the quotient $\C^2/S^1$ is not a manifold. A canonical parameterization of the orbit space is given by the Hilbert map and canonical coordinates are given by the generators of the ring of invariants of the action; see \cite[Thm 5.2.9]{choslaut}. In our case, the ring of invariants is generated by 
\beq
R = |A|^2,\qquad S = |B|^2, \qquad N = A\bar{B}  \label{eq:inv}
\eeq
where, once again, we have set $B = A_\zeta$, with a relation 
\beq
H(R,S,N):=RS - N \bar{N} = 0,\quad R\geq0, \quad S\geq0,\label{eq:rel}
\eeq
hold. The generators of the ring of invariants can be found either through a direct computation or using the general theory of invariants as outlined in \cite[\S5]{choslaut}. Explicitly, one checks that a monomial $A^{\ell_1}\bar{A}^{\ell_2} B^{\ell_3} \bar{B}^{\ell_4}$ is invariant if 
$ \ell_1 + \ell_3 = \ell_2+\ell_4$, which defines a submodule $\mc{M}$ of $\mb{N}^4$. One then readily verifies that this submodule is generated by the vectors $$e_1 = (1,1,0,0), \quad e_2 = (0,0,1,1), \quad e_3 =(1,0,0,1)$$  which correspond to the monomials in \eqref{eq:inv} above. In other words, any invariant monomial can be written as a product of $R,S,N,$ and $\bar{N}$ in an obvious fashion. More generally, one can compute the Molien series \cite[Thm 5.4.1]{choslaut} of this specific representation of the group $S^1$ as $\frac{1+z^2}{(1-z^2)^3}$. The denominator indicates three algebraically independent quadratic monomials, say $R,S,\Re(N)$, and the numerator suggests another, algebraically dependent, quadratic invariant, $\Im(N)$; see \cite[Rem 5.4.2]{choslaut}. 
 

The orbit space $\C^2/S^1$ is then homeomorphic to the subset of $R,S,N\in\R^2\times \C$ where the relations (\ref{eq:rel}) 
hold. One can now express (\ref{s-mfa-n}) in invariants, only, and derive equations for $R,S,N$ \cite[\S 6]{choslaut}. 

We next would like to eliminate the leading-order scaling symmetry, which, since all invariants are quadratic, acts equally on $R,S,N$. Directional blowup \cite{dumortier}, allows us to exploit the leading-order scaling symmetry in an explicit fashion. We therefore use an equivariant, $R$, and associated invariants $S/R$, $N/R$, with respect to the scaling symmetry as new variables. The idea is that the invariants represent the quotient space and the equivariant, which commutes with the action of the scaling group rather than being invariant, tracks the scaling action.  
 One finds that the relations simplify and the orbit space is given as a graph, $S/R=N/R\cdot \overline{N}/R$, so that we may consider the invariant $z=\bar{N}/R=B/A\in\C$, only\footnote{We do not know when to expect such a simplification for more general group actions.}.  We find the system
\begin{align}
z' &= - z ^2 - \hat{c} z - (1 + \ri\hat\omega ) + (1 + \ri \hat\gamma) R\notag\\
R' &= 2 \mathrm{Re}(z R) = R (z + \overline{z})\label{eq:sph}
\end{align}
in the phase space $ \mb{R}^{+}\times\overline{\mb{C}}$. In order to obtain a complete set of charts in a neighborhood of the origin, one also uses the directional blowup in the $S$-direction, with variables $\tilde{z}=A/B=1/z$ and $S$, 
\begin{align}
\tl z' &=1 + \hat{c}\tl z  + \tl z^2\lp[1 + \ri \hat\omega  - (1+\ri \hat\gamma) |\tl z| S^2\rp]\notag\\
S' &= 2S \lp[ -\Re((1-\ri\hat\omega)\tilde{z})-c+\Re((1+\ri\hat\gamma)\tilde{z})|\tilde{z}|^2 S\right)
.\label{eq:psph}
\end{align}
Note that the $\tl z$ equation is precisely the Poincar\'e inversion of the $z$-equation in system \eqref{eq:sph}, and that the singularity $z=\infty$ is non-degenerate, with $\tilde{z}'=1$.

Also, note that the sphere $z\in\C$, $R=0$, together with the point $\tilde{z}=0$, $S=0$, is flow-invariant, the blown-up origin of the original system. On this \emph{singular sphere} $\mc{S}$, we isolated the scaling-invariant, leading-order part of the equation. Of course, the sphere can also be understood as the result of collapsing the Hopf fibration, obtained via $\C^2/S^1\sim \R^+\times (S^3/S^1)\sim \R^+\times S^2$.

Since a direct calculation using \eqref{e:linp} and \eqref{e:scal2} shows that $\hat{c}_\rlin = 2$, it is convenient to introduce the detuning of the trigger speed from the free front speed as a new parameter,
$$
\Delta \hat{c} = \hat{c}_\rlin - \hat{c} = 2 - \hat{c}, \quad0<\Delta\hat{c}\ll 1.
$$
The system \eqref{eq:sph} now reads
\begin{align}
z' &= - (z + 1) ^2 + \Delta \hat{c}\: z - \ri\hat\omega  + (1 + \ri \hat\gamma) R\notag\\
R' &= 2 \mathrm{Re}(z R) = R (z + \overline{z}).\label{eq:sph1}
\end{align}

To conclude this section, we study the flow on the singular sphere, given by the Riccati equation
\begin{equation}\label{e:riccati}
z' = - (z + 1) ^2 + \Delta \hat{c}\: z - \ri\hat\omega  .
\end{equation}
Equilibria 
$$
z_{1/2}^-  = -\fr{2 - \Delta \hat{c}\: }{2}\pm \sqrt{-\Delta \hat{c}\:  - \ri\hat\omega  + \fr{(\Delta\hat{c}\,)^2}{4}},
$$
correspond to eigenspaces of the original linearized equation. The equilibria undergo a complex saddle-node at $z=-1$ when $\Delta \hat{c}\,=0, \,\hat\omega=0$, and a Hopf bifurcation at $\Delta\hat{c}>0$, $\hat\omega=0$. At the Hopf bifurcation, the flow on $\mc{S}$ consists of periodic orbits, whereas outside of the Hopf bifurcation, all non-equilibrium trajectories converge to the same equilibrium. At the saddle-node, all trajectories are homoclinic to $z=-1$. Of course, the Riccati equation can be integrated explicitly, and we will exploit this later on.

We also note that for the $\zeta>0$ dynamics ($\chi(\zeta) = -1$) the equilibria on $\mc{S}$ satisfy
$$
z^2+ (2-\Delta \hat{c}\, ) z + (-1 + \ri\hat\omega ) + \fr{2}{m^2}(1+\ri\alpha), =0.
$$
The equilibrium with negative real part corresponds to the tangent space of the stable manifold $W^{s}_+$ and is explicitly given through
\beq\label{e:z+}
z_+ = -\fr{2 - \Delta \hat{c}\, }{2}- \sqrt{2 - \Delta \hat{c}\, - \rmi\hat\omega + \fr{(\Delta \hat{c}\,)^2}{4}+ \fr{2}{m^2}(1+ \ri\alpha)}.
\eeq

\subsection{Existence of generic free invasion fronts and the blowup geometry}
\label{s:3.4}
We show that generic free fronts exist when $|\alpha-\gamma|\ll 1$ and analyze asymptotics in our blowup coordinates. As we saw in the previous section, the dynamics on the sphere consist of homoclinic orbits converging to a saddle-node equilibrium $z_\rb=-1$ when $\Delta\hat{c}\, = \hat\omega = 0$, corresponding to the unscaled parameters $c=c_\mathrm{lin}$ and $\omega=\omega_\mathrm{lin}$. In particular, there exists a unique trajectory in the strong stable manifold of the equilibrium, which corresponds to decay $\rme^{-\zeta}$. All other trajetories decay with rate $1/\zeta$ in the tangent space of the sphere, which translates into decay $A(\zeta)\sim \mathcal{A}_\infty \zeta\rme^{-\zeta}+\mathcal{B}_\infty\rme^{-\zeta}$ with $\mathcal{A}_\infty\neq0$. In fact, one readily obtains
\beq\label{eq:blowup}
z^f(\zeta) \sim \nu + \fr{1}{\zeta+ \fr{\mc{B}_\infty}{\mc{A}_\infty}} = -1 +\fr{1}{\zeta+ \fr{\mc{B}_\infty}{\mc{A}_\infty}} , \quad R^f(\zeta) = \re^{2\rre\,\nu \zeta} |\mc{A}_\infty \zeta + \mc{B}_\infty|.    
\eeq
The following proposition shows that $\mc{A}_\infty\neq 0$ for $\alpha\sim\gamma$.

\begin{Proposition}\label{prop:genhet}  
For fixed $\alpha$, generic free fronts as defined in Definition \ref{d:1} exist when $|\alpha-\gamma|$ is sufficiently small. 
\end{Proposition}
\begin{Proof}
In scaled variables, when $ \hat{c}=\hat{c}_\rlin = 2$ and $\hat\omega=\hat\omega_\mathrm{lin} =0$, we find 
\beq\label{eq:rgl}
A'' = - A - 2A'+ (1+\rmi \hat\gamma)A |A|^2,
\eeq
where,  from \eqref{e:scal2}, $\hat\gamma=\rmO(|\gamma-\alpha|)$. For $\hat\gamma=0$, this equation posseses a real heteroclinic solution connecting $A=1$ to $A=0$, with the desired asymptotics $A(\zeta)\sim \mathcal{A}_\infty \zeta\rme^{-\zeta}+\mathcal{B}_\infty\rme^{-\zeta}$, $\mathcal{A}_\infty> 0$. Since the heteroclinic is a saddle-sink connection between the circle of relative equilibria $|A|=1$ and the origin, it persists for small values of $\hat\gamma$ as a heteroclinic orbit between a nearby relative equilibrium and the origin. Moreover, the heteroclinic is not contained in the strong stable manifold of the saddle-node equilibrium at $\alpha=\gamma$ and therefore does not lie in the strong stable manifold for $|\alpha-\gamma|\ll 1$.  This proves the proposition.
\end{Proof}

\begin{figure}[h!]
\begin{subfigure}{0.48\textwidth}
\centering
  \def\svgwidth{280pt}
  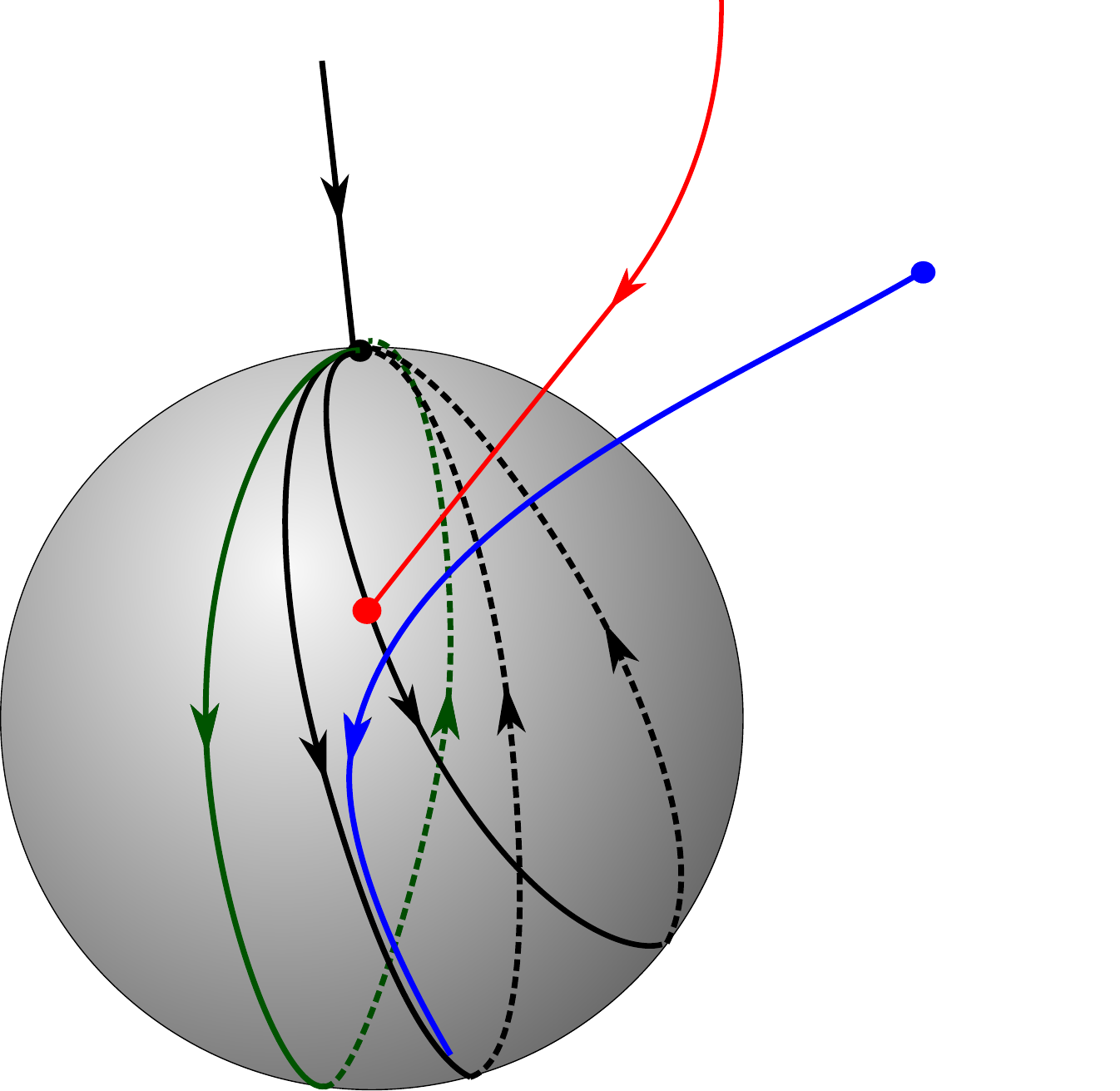
  \end{subfigure}
  \hspace{0.1in}
  \begin{subfigure}{0.48\textwidth}
\centering
  \def\svgwidth{280pt}
  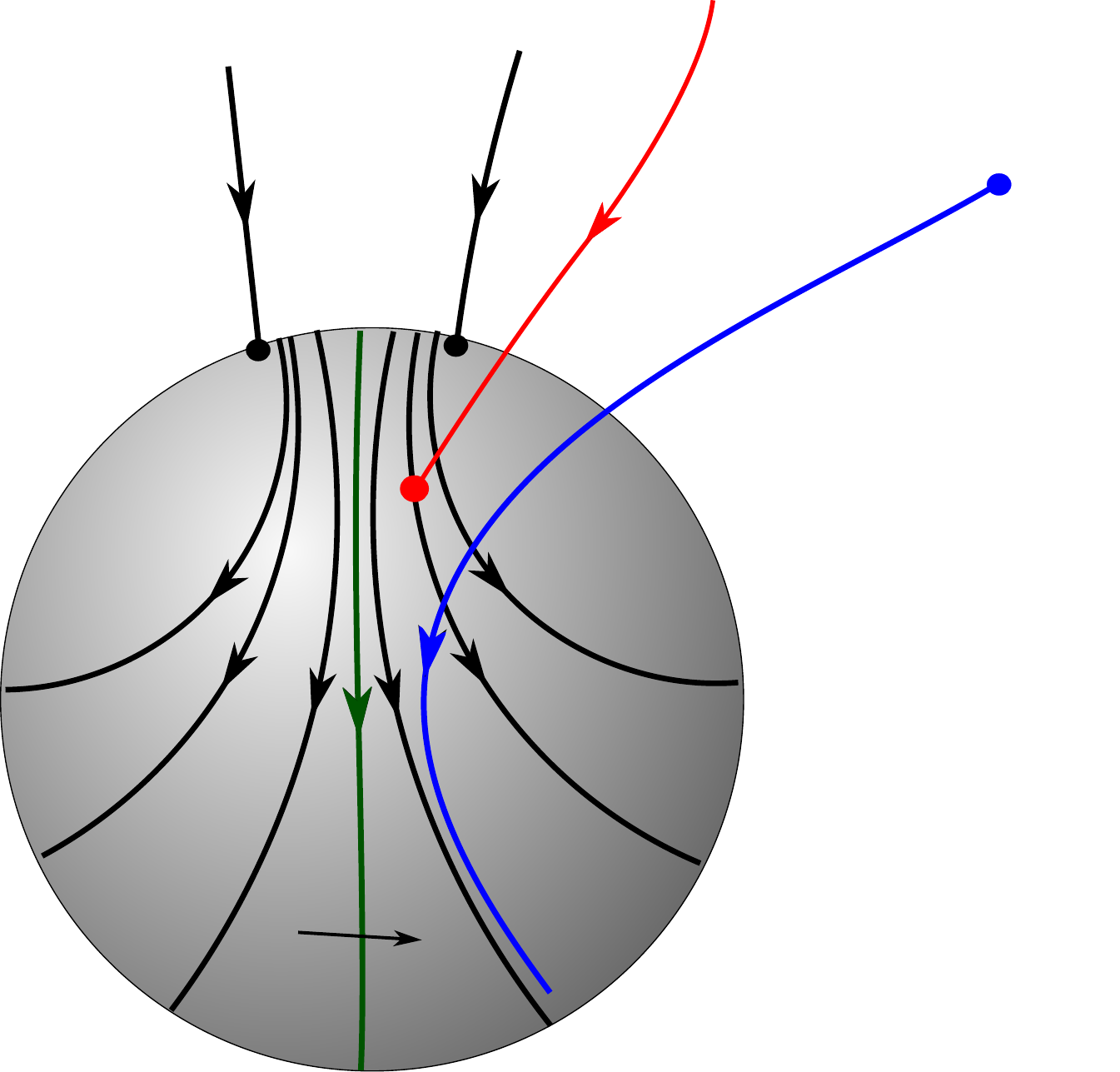
  \end{subfigure}
  \caption{Dynamics near the singular sphere: free front heteroclinic (blue) and  stable manifold $W_+^\rs(0)$ (red). Left: Unscaled parameters set $c = c_\rlin$, $\omega = \omega_\rlin$. On the sphere, homoclinic trajectories are tangent to the real circle (green) at the branch point $z_\rb$.  Right:  Perturb in $c<c_\rlin$ at $\omega = \omega_\rab(c)$.  The branch point $z_\rb$ bifurcates into two equilibria which are encircled by periodic orbits.  For $\omega$ perturbed away from $\omega_\rab$ equilibria become unstable and stable respectively causing drift along the family of periodics.}\label{dia:sph}

\end{figure}

These genericity assumptions can be visualized in the blowup coordinates coordinates of Section \ref{s:3.3}. Since the gauge symmetry is eliminated in such coordinates, the neutral direction of the periodic orbit is removed so that $ \dim W_-^\rcu(A^\rpp) = 1$; see Figure  \ref{dia:sph} for a schematic drawing of dynamics in blowup coordinates near the origin.

\begin{figure}[h!]
\begin{subfigure}{0.48\textwidth}
\includegraphics[width=\textwidth]{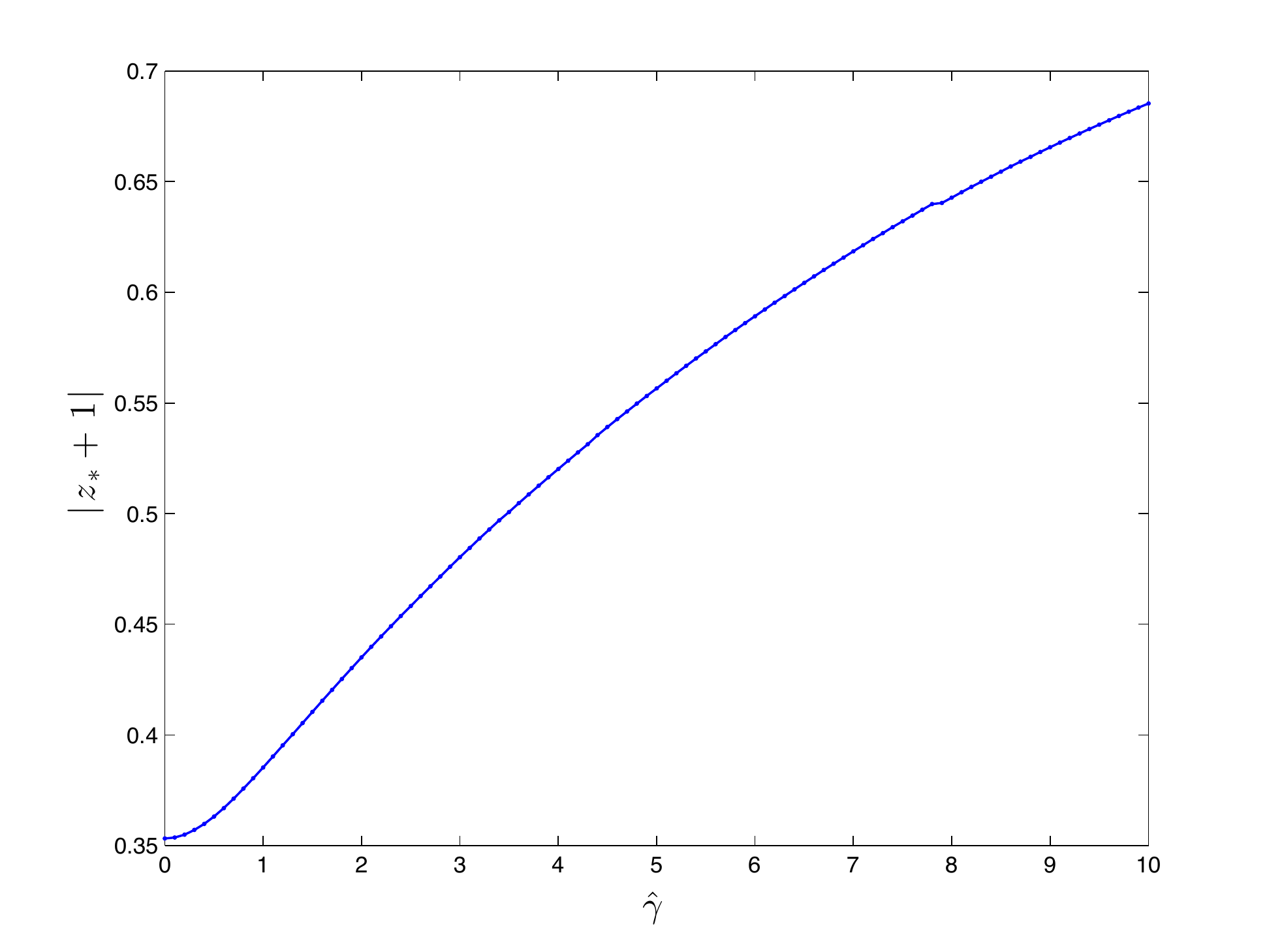}
\end{subfigure}
\begin{subfigure}{0.48\textwidth}
\includegraphics[width=\textwidth]{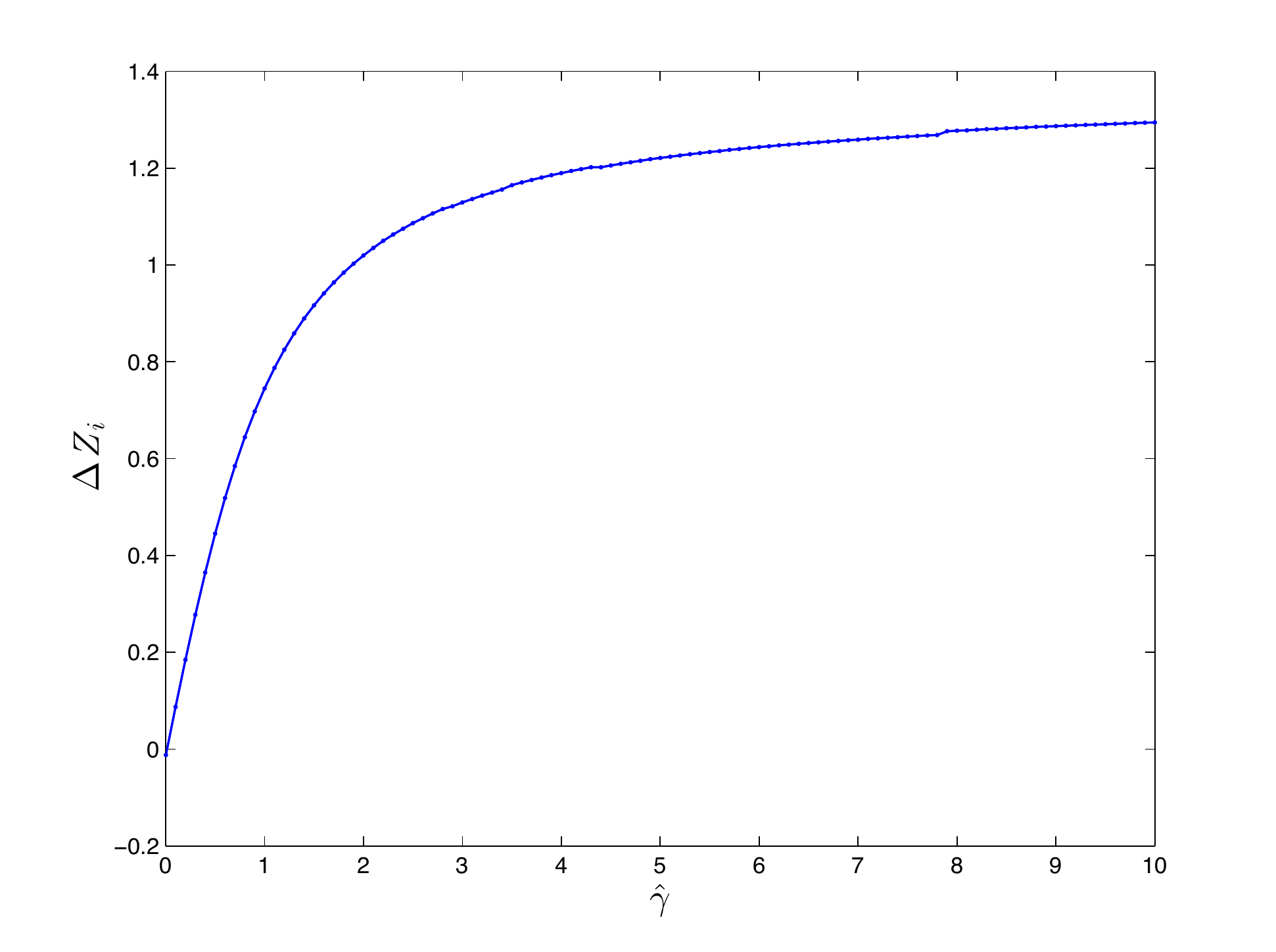}
\end{subfigure}
\caption{Left: Plot of $|z_* +1|$ for $\hat\gamma$ ranging from zero to ten. Right: Plot of $\Delta Z_\ri$ over the same range of $\hat\gamma$.  }
\label{fig:genpo}
\end{figure}

We also tested our hypotheses on existence and genericity of free fronts numerically for $\alpha-\gamma$ not necessarily small. Using a shooting method for \eqref{eq:sph1}, we calculate the unstable manifold of 
$$z_p = \ri \tl k_\rlin,\qquad R= \sqrt{1 - \tl k_\rlin^2},$$ where $\tl k_\rlin$ is the linearly selected wavenumber calculated using \eqref{s-mfl}. \hspace{-0.02in}We then track the base point $z_*$ of the fiber in which this trajectory lies when  $R = \delta$ small.  For parameters at the branch point, if $z_*$ is far away from $z_\rb = -1$, we observe that $W^\rcu_-(A^\rpp)$ does not approach the sphere along the strong stable manifold of the branch point.   Thus if the quantity $|z_*+1|$ is non-zero we obtain that the free front is generic and thus $\mc{A}_\infty\neq 0$.  This quantity is plotted for a range of $\gamma$-values in Figure \ref{fig:genpo}.  Recall that $\alpha = \gamma$ implies that in our scaled coordinates $\hat\gamma = 0$.  Using $z_*$ we can also calculate $\Delta Z_\ri$ defined in (\ref{def:DI}).  We then use this in calculating our $\mc{O}((\Delta c)^{3/2})$-prediction in Figure \ref{fig:autodirect}.  The righthand side shows the coefficient $\Delta Z_\ri$ over a range of $\hat\gamma$-values.

\subsection{Matching stable and unstable manifolds at $\zeta=0$}
\label{s:3.5}

We now prove our main result. The road map is as follows.  We shall see that the singular sphere $\mc{S}$ is normally hyperbolic.  This will give that the flow is locally foliated by smooth one dimensional fibers with base points in $\mc{S}$. Then, well-known results of Fenichel and an inclination lemma will essentially reduce the connection problem to the singular sphere $\mc{S}$. On the sphere, we give a series of coordinate changes which allow us to explicitly integrate the dynamics on $\mc{S}$ in a robust fashion,  dependent on the scaled parameters $\Delta \hat{c}\,$ and $\hat\omega$. We find a connecting orbit using time and $\hat\omega$ as matching parameters. A local analysis yields a leading order expansion for the angular frequency $\hat\omega_{\mathrm{sel}}(\Delta\hat{c})$ in $\Delta\hat{c}$ near zero. We finally unwind our scaling transformations to obtain an expansion of $\omega_\mathrm{tf}$ in the unscaled $ \Delta c\, := c_\rlin -  c$. This expansion can then be inserted 
into the dispersion relation \eqref{eq:nldisp} to obtain a wavenumber prediction for the pattern left in the wake of a trigger front.  

\paragraph{Normal hyperbolicity and smooth foliations near $\mathcal{S}$.} 
In what follows, consider all manifolds in the reduced phase space $\mb{R}_+\times S^2$.
Let $\Phi_\zeta$ be the flow of the full system \eqref{eq:sph1} and $\phi_\zeta$ the flow on $\mc{S}$.
Now, note for $\Delta \hat{c}\: =\hat\omega=0$, all trajectories on the sphere $\mc{S}$ converge to the branch point $z_\rb$ with algebraic decay $\rmO(1/\zeta)$. Since the the linearization at $z_\rb$ is normally hyperbolic, this implies that the invariant singular sphere $\mc{S}$ is normally hyperbolic and Fenichel's results  \cite{Fenichel71,Fenichel77} imply that the phase space $\mb{R}^{+}\times S^2 $ is locally smoothly foliated by smooth one-dimensional strong-stable invariant  fibers $\Phi_\zeta\mc{F}_z\subset \mc{F}_{\phi_\zeta(z)}$. More precisely, all leaves are $C^k$ manifolds, with $C^k$-dependence on base point and parameters, for any finite fixed $k<\infty$. 

As a consequence, there exists a smooth change of coordinates that straightens out the local stable foliation, so that 
\eqref{eq:sph1} can be written as
\begin{align}
z' &= - (z + 1) ^2 + \Delta \hat{c}\:  z -\ri \hat\omega  \notag\\
R' &=R\cdot g(R,z;\hat\omega,\Delta \hat{c}\: ),\label{eq:sphstr}
\end{align}
where $g(0,-1;0,0)=-1$.

To prove existence of the desired connection, we need to choose $\Delta \hat{c}\,>0$ small and find $\hat\omega_*$ and $\zeta_*$ such that there exists a point $(R,z)=(\delta,z_-)\in W^{\rcu}_-(A^\rpp)$ with $\delta>0$ sufficiently small so that 
$$
\Phi_{-\zeta_*}\lp(W_+^\rs(0)\rp) \cap W_-^{\rcu}(A^\rpp)\ni (\delta,z_-).
$$

\begin{figure}[h!]
  \centering
  \def\svgwidth{450pt}
  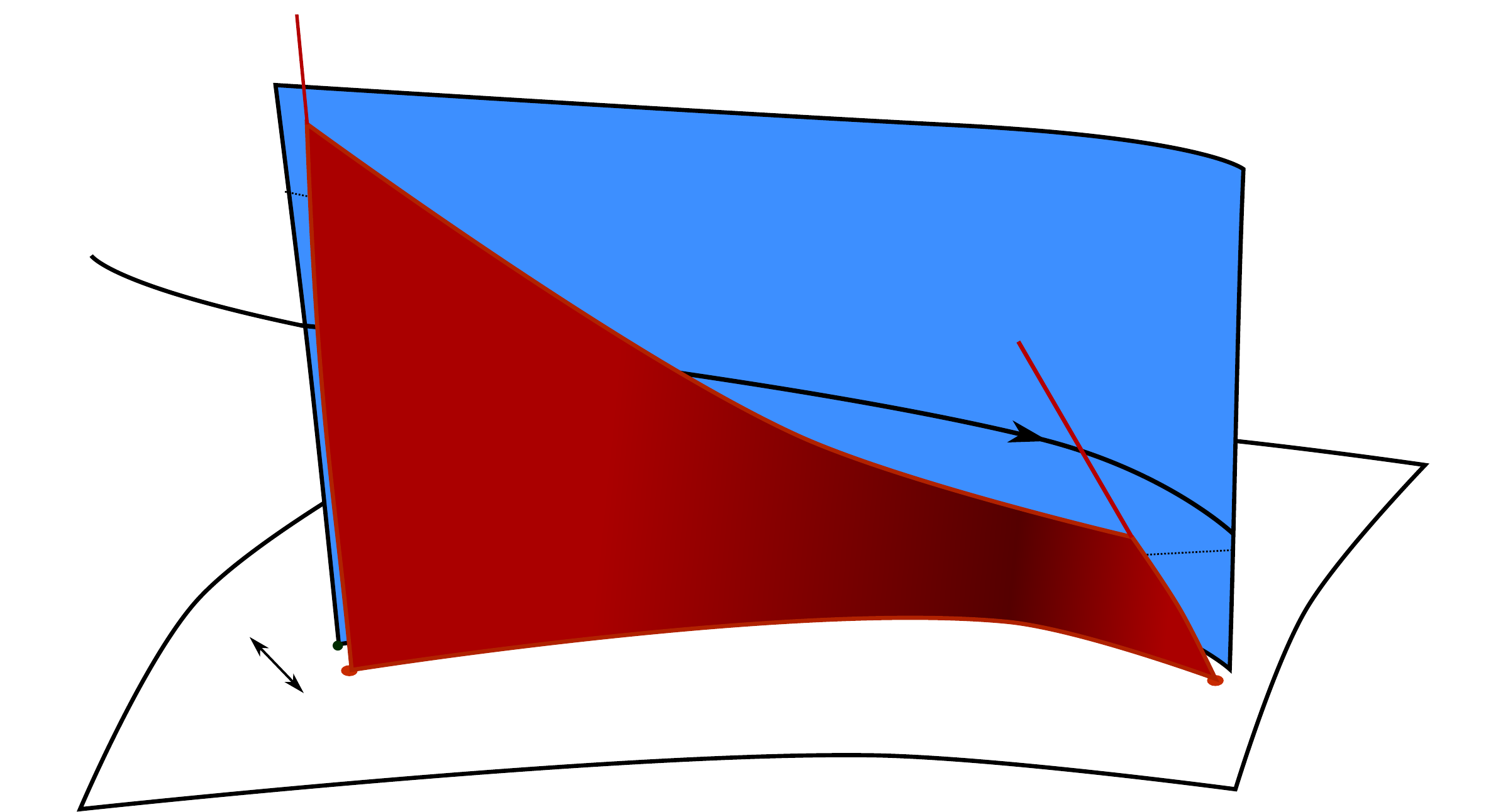
  \caption{Depiction of foliations and heteroclinic connection.  The strong-stable foliation of the trajectory $\phi_\zeta(z_-)$ is denoted as $\mc{M}_{z_-}$ (blue).  Here $\tl h  = \tl h(\delta,-\zeta_*, z_+)$ and $h = \tl h(\delta,0, z_+)$. The stable manifold $W_+^\rs(0)$ (dark red) is stretched in the normal direction under the backwards flow, approaching $\mc{M}_{z_-}$.  We vary $\hat\omega$ and $\zeta$ so that $\Phi_{-\xi_*}(W^\rs_+(0))\cap \mc{F}_{z_-} = (z_-,\delta)$. }\label{dis:full}
\end{figure}

Before we prove the existence of such an intersection,  we investigate the structure of $W_+^\rs(0)$.  It is important to remember that $W_+^\rs(0)$ is not a stable or unstable manifold in the $\zeta<0$ dynamics, but just the set of initial conditions for the $\zeta<0$ flow that will give rise to decaying solutions when integrated forward in the $\zeta>0$ dynamics. In the blowup coordinates $\eqref{eq:sph1}$, $W_+^\rs(0)$ is one dimensional and can be written locally as a graph over the fiber $\mc{F}_{z_+}$ in the $\zeta<0$ dynamics. More precisely, there exists an $h: \mb{R}^{+}\rightarrow \mb{C}$, with $h(0) = 0$ so that locally
$$
W_+^\rs(0) = \{(R, z_+ + h(R))\,|\, R\geq 0\}.
$$
The normal hyperbolicity of $\mc{S}$ then allows us to study how this set evolves under the flow.  Heuristically, when flowed in backwards time, a piece $W_+^\rs(0)_{\mathrm{loc}}$ of the $\zeta>0$ local stable manifold will be stretched out in the normal direction while expanding comparatively little in the $\mc{S}$ directions.  Since $W_+^\rs(0)$ is a graph over the strong stable fiber $\mc{F}_{z_+}$, the considerations above imply that $\Phi_{-\zeta}(W_+^\rs(0))$ will remain a graph over the strong-stable fiber with base point $\phi_\zeta(z_+)$.  This graph will converge to the strong stable fibration as $\zeta$ increases. This strong inclination result is commonly referred to as a $\lambda$-Lemma and is stated precisely in the following lemma; see also Figure \ref{dis:full} for an illustration.

\begin{Lemma}[$\lambda$-Lemma]\label{lem:inc}  The image of the 
 local $\zeta>0$-stable manifold  $W_+^\rs(0)_{\mathrm{loc}}$ under the flow  $\Phi_{-\zeta}$, with $\zeta>0$, large, is exponentially close to the strong stable fiber $\mc{F}_{\phi_{-\zeta}(,z_+)}$.  More precisely, there exist constants $\kappa, C>0$ and a function $\tilde{h}(\cdot,z_+): \mb{R}^{+}\rightarrow \mb{C}$ such that
$$
\Phi_{-\zeta}(W_+^\rs(0)_{\mathrm{loc}}) = \{ \lp(R,\phi_{-\zeta}(z_+) + \tl{h}(R,-\zeta,z_+)\rp)\,\,|\, R\geq0\}
$$
where $\|\tilde{h}(R,-\zeta, z_+)\|_{C^k}\leq C \re^{-\kappa \zeta}$ for any $k<\infty$, $\tl{h}(0,\zeta,z) = 0$, and $\p_z \tl{h}\neq0$. 
\end{Lemma}
\begin{Proof}
This lemma is a consequence of the existence proof for invariant foliations \cite{Fenichel77}: one shows that trial foliations over the linear strong stable foliation converge to the invariant foliation when transported with the backward flow. 
\end{Proof}

�

\paragraph{Connecting via the singular flow.}
With these results in hand,we wish to establish a connection by finding a time-of-flight $\zeta_*$ and frequency $\hat\omega_*$ such that
\begin{equation}
z_- = \phi_{-\zeta_*}(z_+;\hat\omega_*) + \tilde{h}(\delta,-\zeta_*, z_+).\label{eq:con}
\end{equation}
Recall that $z_+$ is the base point of $W^s_+(0)$ and $z_-$ is the base point of the strong stable fiber that contains $W^{\rcu}_-(A^\rpp)$ with $R=\delta$.

We rewrite \eqref{eq:con} as
\begin{equation}
\phi_{-\zeta}(z_+;\hat\omega_*)= z_- - \tilde{h}(\delta,-\zeta^*,z_+),
\end{equation}

and change variables in parameter space,
$$
\mu(\Delta \hat{c}\: ,\hat\omega) = -\Delta \hat{c}\, - \ri\hat\omega + \fr{(\Delta \hat{c}\, )^2}{4}:= M^2 \re^{\ri \theta},\qquad M,\theta\in \mb{R}.
$$
Next, we shift $z = \hat{z}  +  \fr{\Delta \hat{c}\,}{2} - 1$ in \eqref{eq:sph} 
and obtain
\beq
 \hat{z} ' = - \hat{z}^2 + \mu(\hat\omega,\Delta \hat{c}\, ).
\eeq
In these variables, $W_+^\rs(0)$ intersects the blowup sphere at the point
$$
\hat{z}_+=\hat{z} _+(\Delta\hat{c},\hat{\omega})= z_+  - \fr{\Delta \hat{c}\,}{2} + 1 = - \sqrt{K- \Delta \hat{c}\, +\fr{(\Delta \hat{c}\,)^2}{4}- i \hat\omega } =  - \sqrt{K + M^2 \re^{i\theta}},
$$
where $K = 2+ \fr{2}{m^2}(1+\ri\alpha )$, $m$ was defined in (\ref{e:scal}), and $z_+$ was defined in (\ref{e:z+}). 

We note that $\hat{z}_+(0,0) =  4+2\ri\alpha$ is non-zero.  Also, our assumptions give that $\hat{z}_-$ is non-zero, but needs to be evaluated numerically.  This allows us to define, 
\beq\label{def:DI}
\Delta Z = \lp(\fr{1}{\hat{z}_+} - \fr{1}{\hat{z}_-}\rp)\Bigg|_{M=0},\qquad \Delta Z_\mathrm{r} = \rre\,\Delta Z, \qquad \Delta Z_\ri = \rim\,\Delta Z.
\eeq
Note that $M=0$ corresponds to $\Delta\hat{c} = \hat{\omega} = 0$.

Next, we scale $\hat{z} = Mw,\, \hat{\zeta} = M\zeta$ and use the M\"obius transformation 
$$
\rho= \fr{w + \eta}{w- \eta},\qquad \eta = \re^{\ri\fr{\theta}{2}}
$$
to shift equilibria to the north and south pole of the sphere, respectively.

Last, we set $ r = \log(\rho)$ and find the constant vector field
\beq
\dot{r} = (\log(\rho))_{\hat{\zeta}} = 2\eta.\label{eq:const}
\eeq

In our new coordinates, the points $\hat{z}_{+}$ and $\hat{z}_- - \tilde{h}(\delta,\hat{z}_-,-\hat{\zeta}/M)$, which we wish to connect, have the representations
$$
r_{+} =  \log\lp(\fr{1+ \fr{M \eta}{\hat{z}_{+}}}{1 - \fr{M \eta}{\hat{z}_{+}}} \rp) + \ri(2\pi j_{+}) , \qquad r_{-} =  \log\lp(\fr{1+ \fr{M \eta- \tl h}{\hat{z}_{-}}}{1 - \fr{M \eta-\tl h}{\hat{z}_{-}}} \rp) + \ri(2\pi j_{-}), 
$$
where $\tl h = \tl h(\delta,\hat{z}_+, -\hat{\zeta}/M)$, and $j_{\pm}\in \mb{Z}$ take into account that the complex logarithm is multi-valued.
By varying $\mu$ (i.e. $\hat\omega$ and $\Delta \hat{c}$) we wish to find a solution $r(\widehat{\zeta}\,)$ connecting these points in finite time $\widehat{\zeta} = \widehat{T}$.  

As we are expanding from $\Delta \hat{c}\:  = \omega = 0$, we have that $M$ is small. Thus we find to first order
$$
r_{+} = \ri(2\pi j_{+})+ \fr{2M \eta}{\hat{z}_{+}} + O(M^{3}),\qquad
r_-  = \ri(2\pi j_-) + 2\fr{M\eta - \tl h}{\hat{z}_-} + \rmO(M^3).
$$

Now integrating \eqref{eq:const}, setting $r(0) = r_+$ and $r(-\hat{\zeta}) = r_-$,  we obtain
\beq
-2\eta\, \hat{\zeta}= r_- - r_+ =- \ri(2\pi \Delta j) - 2 \eta \,M\,\Delta Z  -  \fr{2 \tl h}{\hat{z}_-} +  O(M^3),
\eeq
where $\Delta j = j_+ - j_-$ and we have discarded all but the leading order term, in $M$, of $\fr{1}{\hat{z}_+} - \fr{1}{\hat{z}_-}$.  We then obtain the equation
\beq
\hat{\zeta} =  \ri\Delta j\pi  \re^{-i\fr{\theta}{2}} + M \cdot \Delta Z+ \fr{2 \tl h}{\eta\hat{z}_-} + O(M^3).\label{eq:time}
\eeq
By Lemma \ref{lem:inc}, we have for $\hat{\zeta}>0$, $\tl h = \tl h(\delta, -\hat{\zeta}/M,\hat{z}_+) = \rmO(M^k)$ for all $k\geq1$.  Thus, we can smoothly extend $\tl h$ at $M = 0$ to $\tl h = 0$. 

Now setting $\hat{\theta} = -\pi + \theta$ we obtain
\beq
\hat{\zeta} =  - \Delta j\pi  \re^{-\fr{\hat{\theta}}{2}} + M \Delta Z+ O(M^3)\label{eq:sphcon}
\eeq
This equation has the solution
$$
M = 0, \quad \hat{\theta} = 0,\quad \hat{\zeta} = -\Delta j \, \pi.
$$
We choose $\Delta j<0$ so that $\hat{\zeta}>0$. We wish to find a solution near this point for $M$ small.   Considering the imaginary part of \eqref{eq:sphcon}
$$
0 = \Delta j \pi \sin (\hat{\theta}/2) + M\,\Delta Z_\ri + \rmO(M^3)
$$
and expanding in $\hat{\theta}$, we obtain
$$
\hat{\theta} = -2 \fr{M\,\Delta Z_\ri}{\pi \Delta j} + \rmO(M^3).
$$
Inserting this into the real part of \eqref{eq:sphcon} and solving for $\hat{\zeta}$, we obtain
$$
\hat{\zeta} = -\pi \Delta j+ M\,\Delta Z_\mathrm{r} + \rmO(M^2).
$$

Now using the polar change of coordinates $-\Delta \hat{c}\, + (\Delta \hat{c}\, )^2/4 - \ri \hat\omega = M^2 \re^{\ri\theta}$, the Implicit Function Theorem allows us to obtain the expansion 
$$
-\Delta \hat{c}\:  + \fr{(\Delta \hat{c}\, )^2}{4} - \ri \hat\omega = -M^2\cdot\lp(1 + \ri \fr{2\Delta Z_\ri}{\pi\Delta j}M + \rmO(M^2)\rp).
$$
Noticing that $ M = \sqrt{\Delta \hat{c}\,} + \rmO(\Delta \hat{c}\,)$, we solve the imaginary part of this equation for $\hat\omega$ and expand

$$ \hat\omega(\Delta \hat{c}\,) = \fr{2\Delta Z_\ri}{\pi \Delta j} (\Delta \hat{c}\, )^{3/2} + \rmO((\Delta \hat{c}\, )^2).\label{eq:expan}
$$ 
We summarize the above discussion in the following proposition. 
\begin{Proposition}\label{p:exp}
Set $\Delta Z_\ri := \rim(\fr{1}{\hat{z}_+} - \fr{1}{\hat{z}_-})$ where $z = \hat{z} + \fr{\Delta \hat{c}\,}{2} - 1$ and choose $\Delta j\in\Z_-$, fixed.  Then for $\Delta \hat{c}\,>0$ small, there exists $\zeta_*$ and $\hat\omega_*$ for which  \eqref{eq:con} is satisfied.  Furthermore, we have the following expansions:
\beq\label{e:expp}
\hat\omega_*(\Delta \hat{c}\,) = \fr{2\Delta Z_\ri}{\pi \Delta j} (\Delta \hat{c}\, )^{3/2} + O((\Delta \hat{c}\, )^2),\qquad  {\zeta}_* = -\pi (\Delta j) (\Delta \hat{c}\,)^{-1/2}+ \Delta Z_\mathrm{r}   + \rmO((\Delta \hat{c}\,)^{1/2}).
\eeq
\end{Proposition}

\begin{Remark}
We emphasize that the behavior of $W_+^\rs(0)$ and $W_-^\rcu(A^\rpp)$ near the origin determine the coefficient $\Delta Z_\ri$ of the leading order term in the expansion (\ref{e:expp}). As pointed out in the introduction, $\Delta Z_\ri$ measures a distance between the leading edge of the front and the stable subspace in projective coordinates. This can be made explicit by considering the Poincar\'e-inverted coordinates $z_\mathrm{inv}=1/(z+1)$, in which $\Delta Z=z_\mathrm{inv}^+-z_\mathrm{inv}^-$. Since $z_\mathrm{inv}^-$, the base point of the stable fiber corresponding to $W^{cu}_-$, is only defined up to flow translates, we can vary $\Re(z_\mathrm{inv}^-)$ arbitrarily by shifting the free front. Therefore, the imaginary part $\Delta Z_\ri$ merely measures the distance between 
$z_\mathrm{inv}^+$ and the homoclinic orbit in the base of the fibers corresponding to $W^\mathrm{cu}_-$.

Also notice that, since trajectories approaching $\mc{S}$ have the asymptotic form $\eqref{eq:blowup}$,  we have that 
$$
\rim\{\fr{1}{z_-}\} = \rim\{\zeta + \mc{B}_\infty/\mc{A}_\infty\} = \rim\{\mc{B}_\infty/\mc{A}_\infty\}.
$$ 
\end{Remark}

\subsection{Proof of Theorem \ref{t:1}}

Proposition \ref{p:exp}  gives the existence of trigger fronts:   Given $(\hat{\omega}_*, \hat{\zeta}_*)$ there are corresponding $ \omega_*$ and $\xi_*$ such that $W_+^\rs(0)$ and $W_-^{\rcu}(A^\rpp)$ intersect non-trivially.  By taking points in this intersection as an initial condition at $\xi = 0$, then sending $\xi \rightarrow -\infty$, the trajectory (which is in $W_-^{\rcu}(A^\rpp)$) must converge to $A^\rpp$ .  In the same way, sending $\xi\rightarrow +\infty$ the trajectory (which is in $W_+^\rs(0)$) must converge to $A\equiv 0$.

In order to prove our main result, it remains to track the scalings from Section \ref{s:3.1}. We therefore reintroduce hats, writing for instance  $\hat{\omega}$ for the frequency in scaled coordinates as given in Proposition \ref{p:exp}, and $\omega$ for the parameter in CGL \eqref{e:cglt}; see (\ref{e:scal}) and (\ref{e:scal2}). 

Using  \eqref{e:scal2} to write the leading order expansion from the above proposition in terms of the unscaled variables  $\Delta c := ( c_\rlin -  c)$ and $\Delta \omega = (\omega  - \omega_\rab(\alpha, c))$, we obtain
\begin{align}
\fr{\Delta \omega}{m^2} &=
  \fr{2 \Delta Z_\ri}{\pi \Delta j} \fr{(\Delta  c)^{3/2}}{(m^2 (1+\alpha^2))^{3/4}} + \rmO((\Delta c)^2).\notag
\end{align}

Further simplifying we find
\begin{align}
 \fr{2 \Delta Z_\ri}{\pi \Delta j} (\Delta c)^{3/2} +\rmO((\Delta  c)^2)&= \fr{(1+\alpha^2)^{3/4}\Delta\omega}{m^{1/4}}\notag\\
 &= (1+\alpha^2)^{3/4}\Delta\omega + \rmO(\Delta \omega \Delta  c, (\Delta \omega)^2, (\Delta  c)^2),\label{eq:expan1}
 \end{align}
 where we have used the fact that $m^{-1/4} =  1 + \rmO(\Delta c, \Delta \omega).$

From here we can use the Implicit Function Theorem to obtain
$$
\Delta \omega =  \fr{2 \Delta Z_\ri}{\pi \Delta j(1+\alpha^2)^{3/4}} (\Delta  c)^{3/2}  + \rmO((\Delta  c)^{2})
$$
from which we obtain the expansion as stated in the main theorem by setting $\Delta j = -1$:
$$
\omega =  \omega_\rab(\alpha, c) - \fr{2 \Delta Z_\ri(c_\rlin - c)^{3/2}}{\pi \Delta j(1+\alpha^2)^{3/4}}  + \rmO( ( c_\rlin -  c)^{2}).
$$
In the same manner we can also obtain the expansions for the unscaled $\xi_*$, the distance between the trigger and the invasion front,  and $k_\rtf$ the selected wavenumber as given in Theorem $\ref{t:1}$.

\section{Comparison with direct simulations and heteroclinic continuation}\label{sec:num}

\begin{figure}[h!]
\begin{subfigure}{ 0.48\textwidth}
\centering
\includegraphics[width=\textwidth]{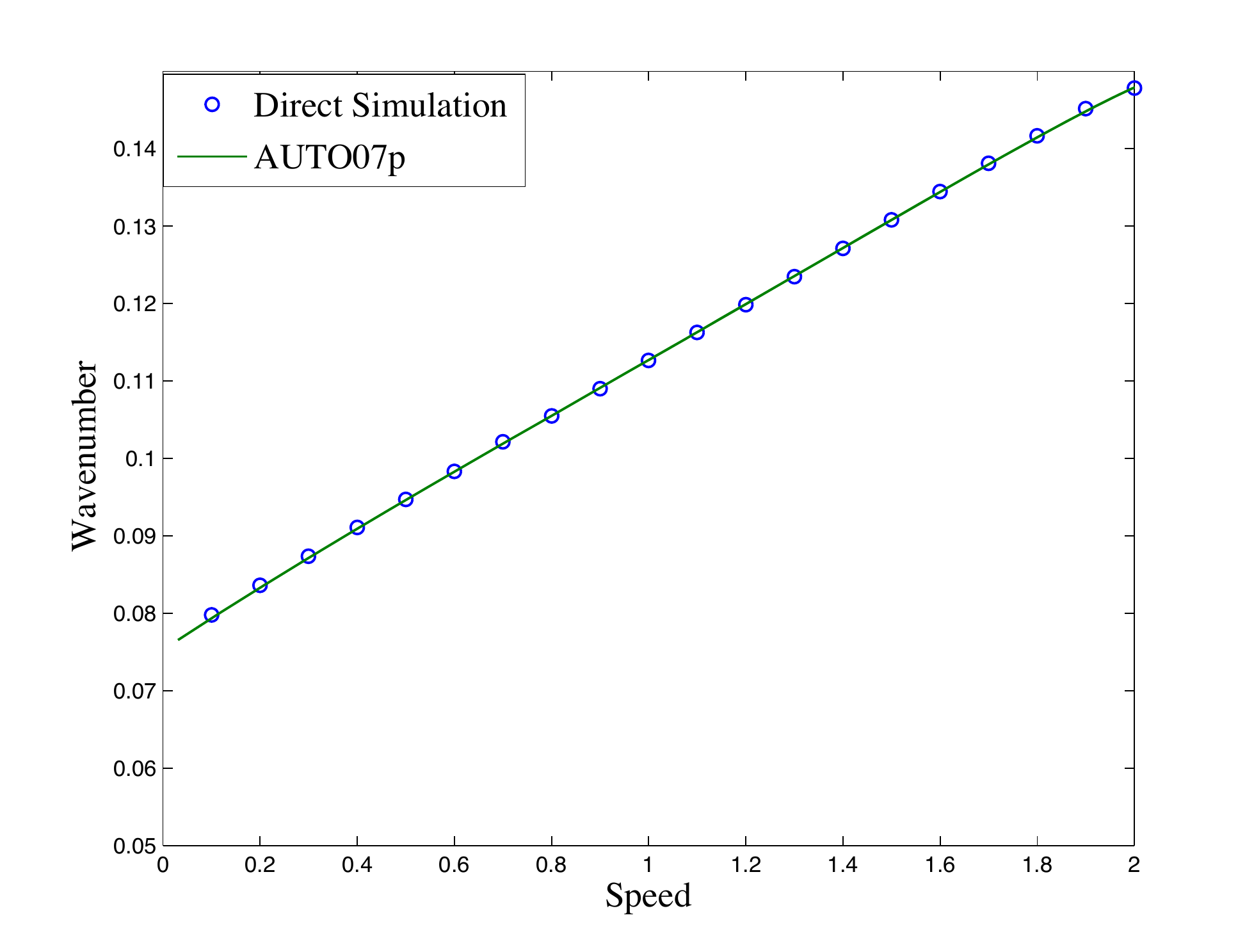}
\end{subfigure}
\begin{subfigure}[h!]{0.48\textwidth}
\includegraphics[width=\textwidth]{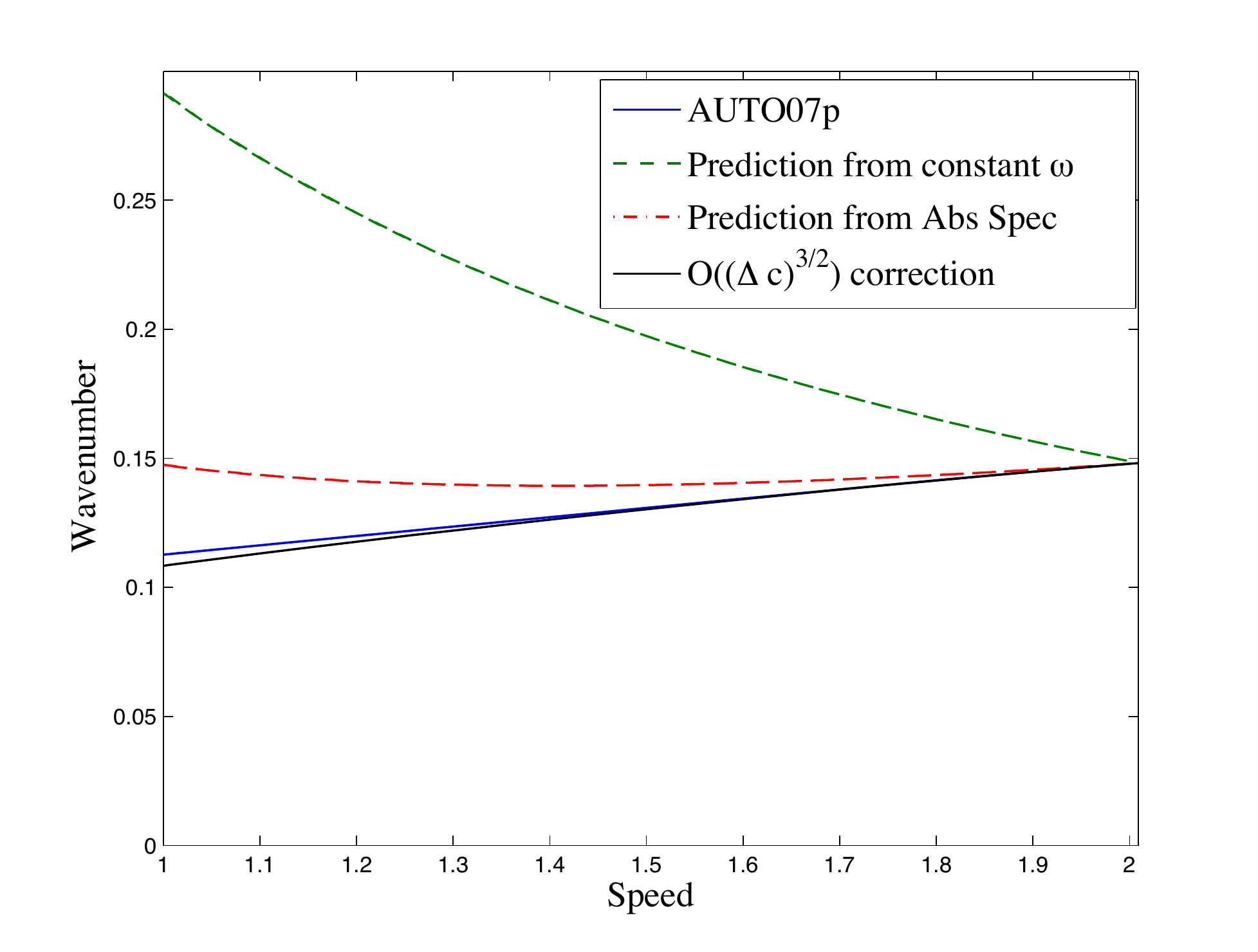}
\end{subfigure}
\centering
\vspace{-0.1in}
\caption{Left: Comparison of wavenumbers from direct simulations with \textsc{auto07p} calculations for a range of $c$ values for fixed $\alpha = -0.1, \gamma = -0.2$. Discrepancy between the two calculations is less than $0.1\%$ for $dt = 0.01$ and $dx = 0.0767$.  Right: Comparison of different predictions for the selected wavenumber with AUTO07p calculations. The speed is varied while other parameters are fixed at $\alpha = -0.1,\gamma = -0.2$.  
}
\label{fig:autodirect}
\end{figure}

We now compare our leading order prediction with both direct simulation of \eqref{e:cglt} and computation of trigger fronts via continuation.  For the direct simulations, we used a first order spectral solver, simulating in a moving frame of speed $c$. The calculations were performed on a large domain ($ L = 2400$) to avoid wavenumber/wavelength measurement error.  Since wavelengths converge slowly as expected \cite{vanSaarloos03}, simulations were run for sufficiently long times ($t \sim 5000$).    
To corroborate these simulations, we used the continuation software \textsc{auto07p} to find the heteroclinic orbit directly as solution to a truncated boundary value problem, and then continue in $c$ with $\omega$ as a free parameter. As seen in Figure \ref{fig:autodirect}, both methods are in reasonably good agreement.

\begin{figure}[h!]
\centering
\begin{subfigure}[h!]{0.5\textwidth}
\includegraphics[width=\textwidth]{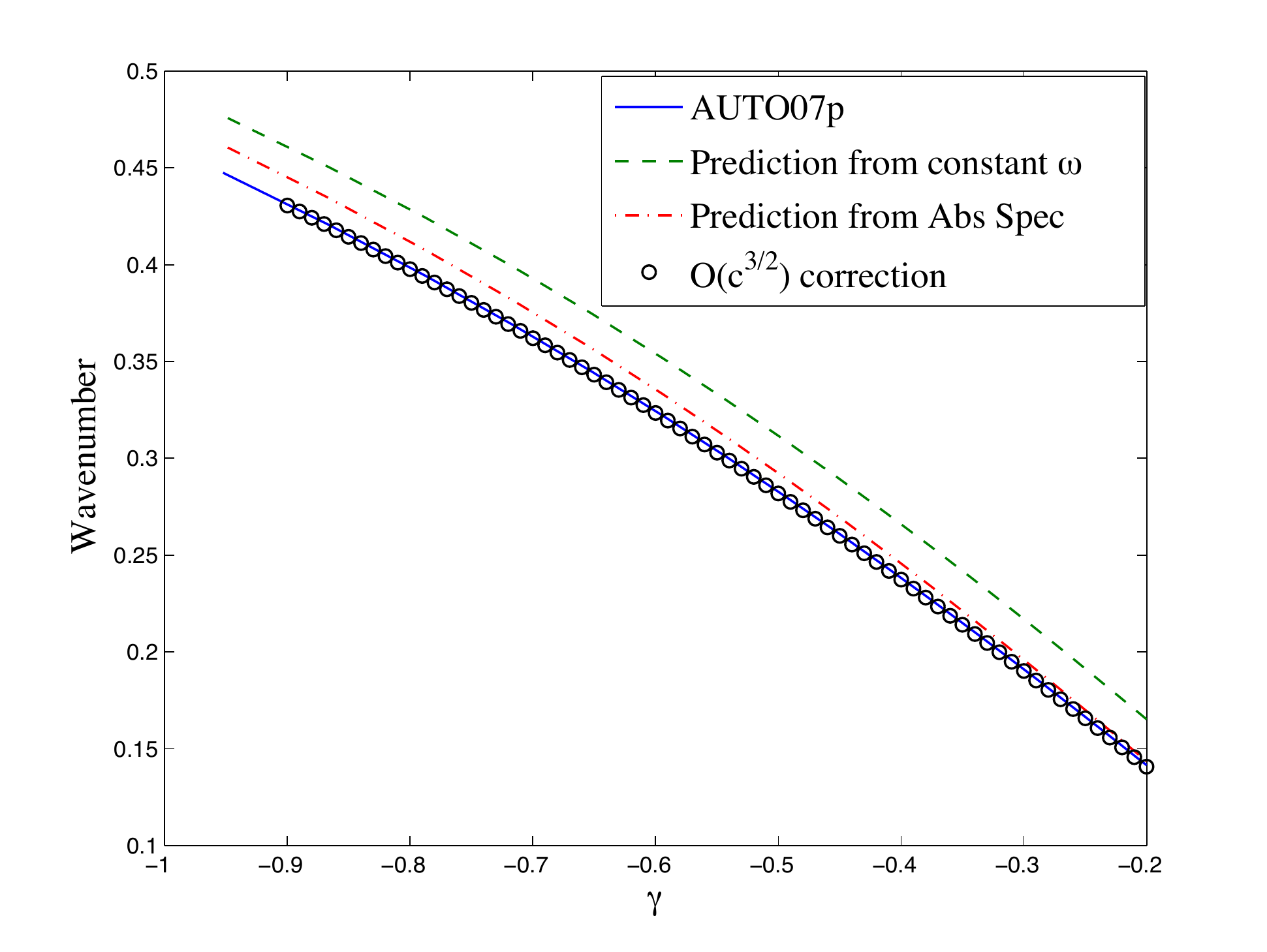}
\end{subfigure}
\hspace{-0.3in}
\begin{subfigure}[h!]{0.5\textwidth}
\includegraphics[width=\textwidth]{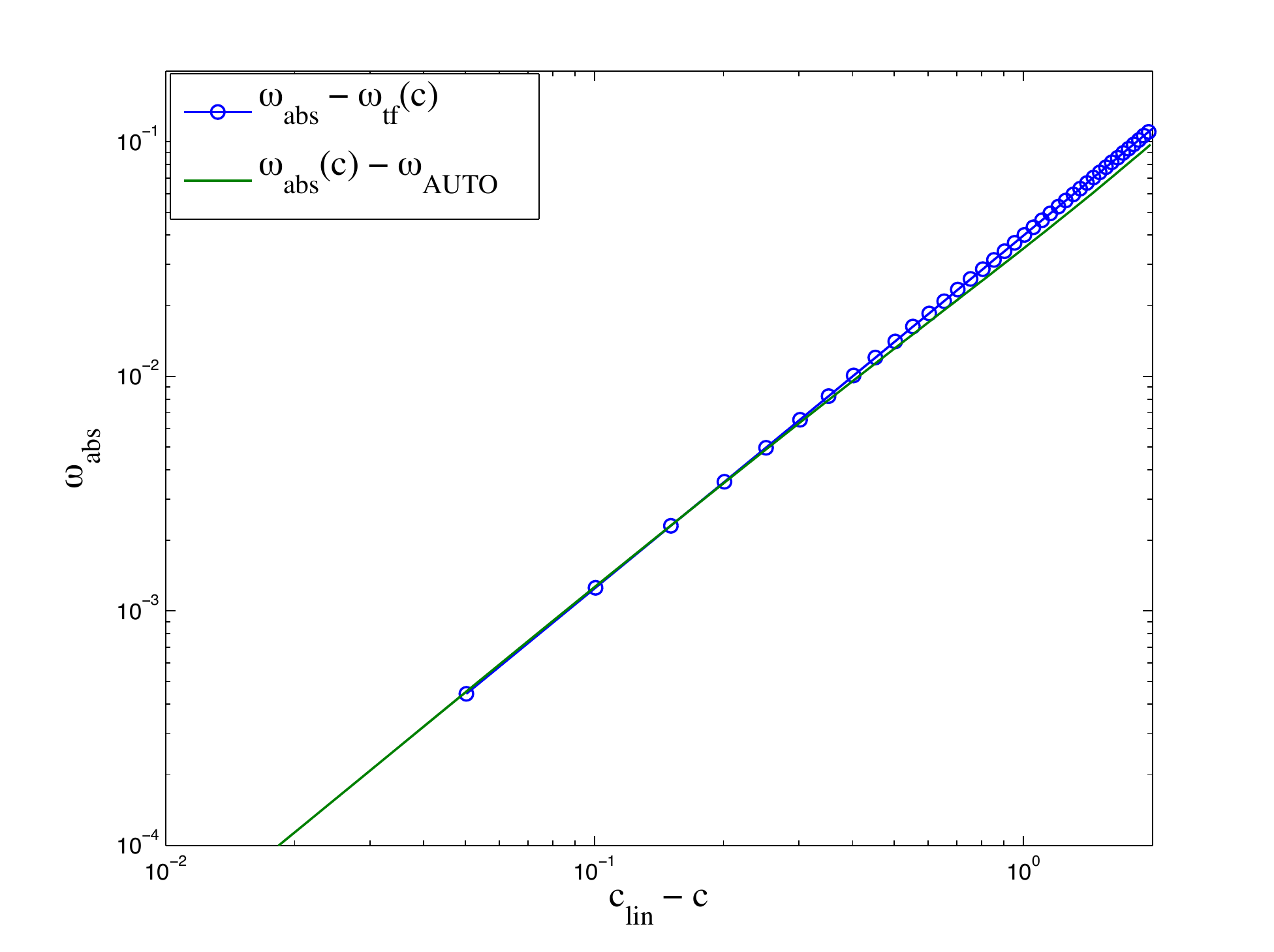}
\end{subfigure}
\vspace{-0.1in}
\caption{Left: Comparison of different predictions for the selected wavenumber with AUTO07p calculations for a range of $\gamma$ with fixed $ c = 1.8$ and $\alpha = -0.1$. Right:  Logarithmic plot of the difference $\omega_\rab(c) - \omega_\rtf(c)$ using expansion from the Theorem \ref{t:1} above, and data from AUTO07p continuation.   }\label{fig:compall}
\end{figure}

We also give comparisons between \textsc{auto07p} computations, the prediction based on the absolute spectrum $\omega\sim\omega_\mathrm{abs}$, and the $\rmO((\Delta c)^{3/2})$ correction from Theorem \ref{t:1} in Figure \ref{fig:autodirect}.  There we also compare with a somewhat naive prediction, assuming that the frequency $\omega$ of the invasion process is constant at leading order. As pointed out in Figure \ref{fig:absspec}, this neglects frequency detuning based on wavenumber changes (slope of absolute spectrum) and speed (shift of absolute spectrum).  Finally, Figure \ref{fig:compall} shows how our predictions fare when $c$ is fixed and $\gamma$ is varied, as well as a log-log plot of $\omega_\rtf - \omega_\rab$, confirming the exponent $3/2$ and the coefficient $\Delta Z_\ri$.

\section{Discussion and Future Work}\label{sec:disc}
 
In this work, we have proved the existence of a coherent triggered front in the CGL equation for trigger speeds close to, but less than, the linear invasion speed, under mild generic assumptions on free invasion fronts.  Furthermore, we have shown how the trigger \emph{selects} the periodic wave-train created in the wake and how its speed affects the wavenumber. Our main tools were a sequence of coordinate changes in a neighborhood of the origin, based on geometric desingularization and invariant foliations. As a byproduct, we establish  expansions for the frequency of the trigger front with universal leading order coefficient determined by the absolute spectrum. At higher order, coefficients in the expansion depend on a projective distance between the leading edge of the front and a stable manifold ahead of the trigger. 
 
While some of the tools used in this paper may not immediately apply in other situations (such as those mentioned in the introduction), we expect that one could adapt the main concepts from Section \ref{sec:dynbp}.  Namely, for a triggered front in such systems, selected wave numbers should be determined by the intersection of the absolute spectrum with the imaginary axis at leading order.  In many of these systems, explicit expressions for the spectra are not known but simple algebraic continuation usually allows one to easily obtain accurate predictions.

The present paper addresses existence and qualitative properties of trigger fronts in the simplest possible (yet interesting) context, leaving many open questions. 

First, it would be interesting to study stability of the trigger fronts. Stability is determined, at first approximation, by spectral properties of the linearization at such a front.  Essential and absolute stability of this linearized operator are determined by essential and absolute spectra at $\xi = \pm \infty$.  Since we have stability at $\xi = +\infty$, the only destabilizing influence in the far field comes from instabilities of wave trains in the wake.  For $ c \sim c_\rlin$ these instabilities are known to be convective  \cite{nozakibekki,vanSaarloos92} or absent for $\alpha,\gamma$ not too large.  In addition, it would be interesting to study and possibly exclude instabilities via the extended point spectrum.  Based on the geometric construction of coherent trigger fronts and direct simulations, we do not anticipate such instabilities.  We do suspect however that fronts with $\Delta j <-1$ (\ref{eq:omselj})  would pick up unstable eigenvalues near $A = 0$ and hence be unstable. A larger $\Delta j$ will lead to fronts with larger distance $\xi_*$ to the trigger location $\xi=0$ (\ref{e:expp}) and several small oscillations in this gap. One then expects this unstable plateau to generate unstable eigenvalues in the linearization; see \cite{gluing} for a similar scenario. From a different perspective, fronts with higher $|\Delta j|$ arise through bifurcations from an already unstable primary state as explained in the discussion of the role of absolute spectra in Section \ref{sec:dynbp}, when considering large, bounded domains. As a consequence, one expects the bifurcated states to be unstable as well. 
   \begin{figure}[h!]
\begin{subfigure}[h!]{0.5\textwidth}
\includegraphics[width=\textwidth]{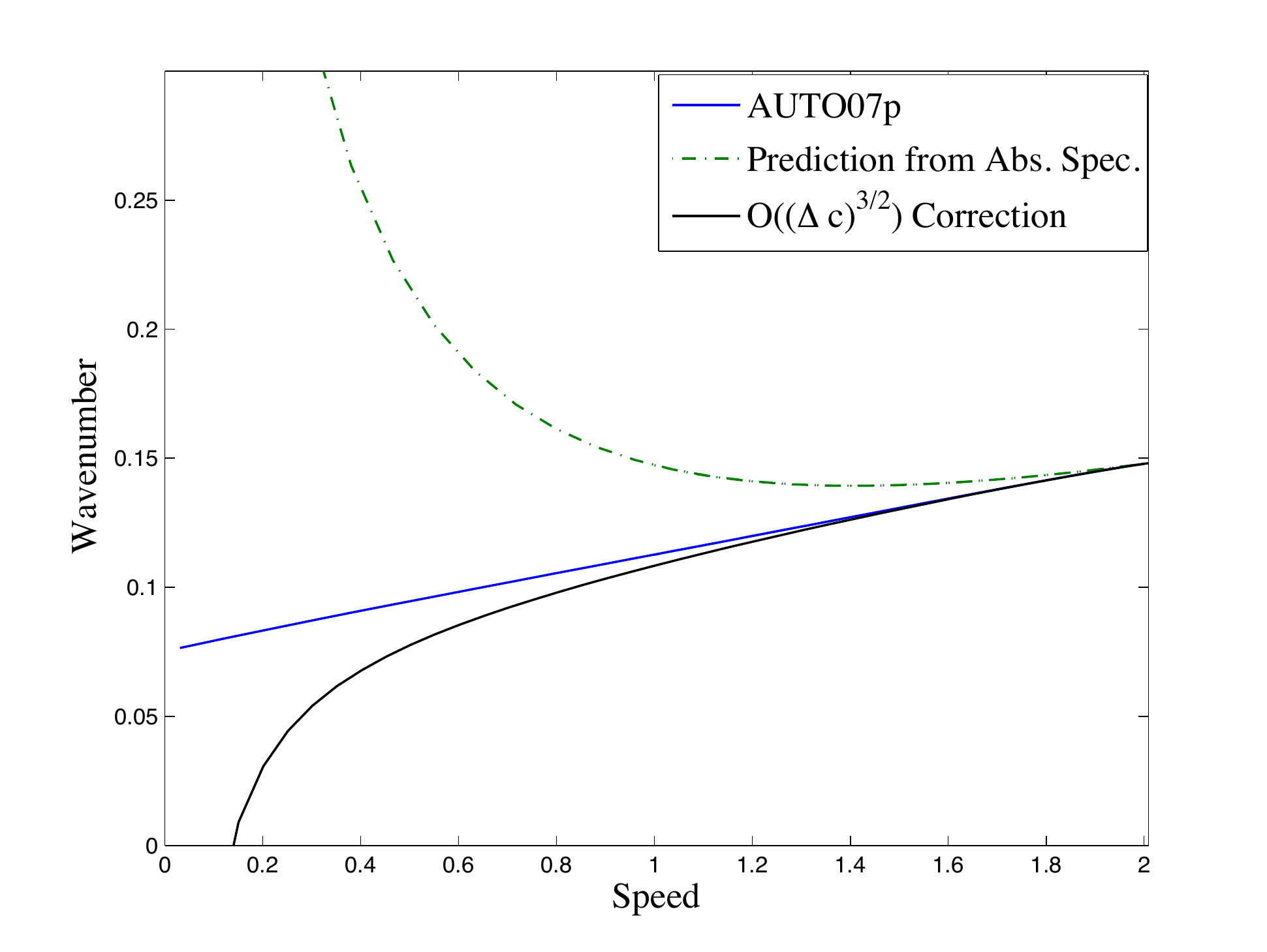}
\end{subfigure}
\begin{subfigure}[h!]{0.5\textwidth}
\includegraphics[width=\textwidth]{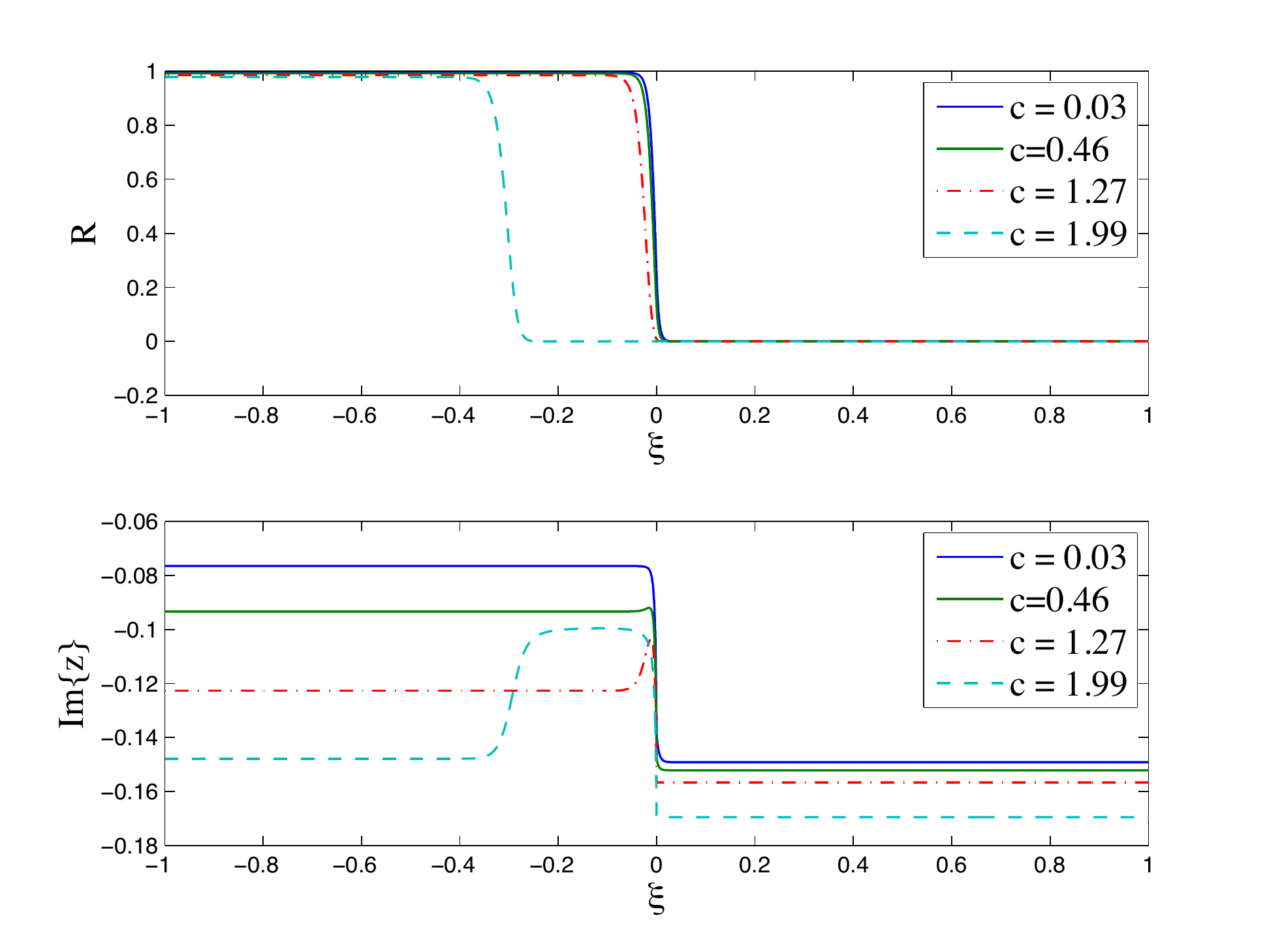}
\end{subfigure}
\caption{We fix $\alpha = -0.1, \gamma = -0.2$. Left: Comparison of AUTO07p calculation with the absolute spectrum prediction and the $\rmO((\Delta c)^{3/2})$ correction for speeds ranging from $c = 0 $ to $ c = c_\rlin$. Right: Front profiles calculated in the blow-up coordinates using AUTO07p for various speeds.  Here the left boundary condition is given by the periodic orbit $(z,R) = (\ri k_\rtf,1 - k_\rtf^2)$ and the right is given by $(z_+,0)$.    }
\label{fig:compallsp}
\end{figure}

One can also envision many generalizations of our result. Using ill-posed spatial dynamics on time-periodic functions as in \cite{Goh11,Scheel03}, one can study similar problems in pattern-forming systems without gauge symmetry. In a different direction,  we expect that triggers $\chi$ with $\chi'$ sufficiently localized would be immediately amenable to our analysis. In particular, monotone triggers with $\chi'$ sufficiently localized should yield the same type of expansion, albeit with different, non-explicit, projective distances $\Delta Z$. On the other hand, one can envision how non-monotone triggers with long plateaus where $\chi>\chi_-$, would generate triggered fronts even for speeds $c>c_\mathrm{lin}$. In terms of our linear heuristics, the linearized operator at the origin, $(1+\rmi\alpha)\partial_{\xi\xi}+c\partial_\xi + \chi-\rmi\omega$ now possesses unstable extended point spectrum \cite{SansSchRad07} in addition to the unstable absolute spectrum. In a different direction, slowly varying triggers $|\chi'|\ll 1$ or triggers that do not simply modify the linear driving coefficient pose a variety of interesting challenges. 

Another direction is suggested by Figures \ref{fig:compall} and \ref{fig:compallsp}.  For speeds further away from the linear spreading speed, our predictions deviate significantly from the actual selected wavenumber and it is not clear in which context one might be able to establish analytic predictions. Going all the way to trigger speed $c = 0$ could however serve as another point in parameter space where analytic expansions can be derived. In fact, thinking of the trigger as an effective Dirichlet-type boundary condition,  one expects standing triggers similar to Nozaki-Bekki holes \cite{standing}. Such holes are explicitly known coherent structures in CGL that emit wave trains.  Selected wavenumbers for $ c$ close to zero would be corrections to the wave numbers selected by these coherent structures.

\bibliography{cgltriggerJNLSedits}

\def\cprime{$'$}
\begin{thebibliography}{10}

\bibitem{mimura}
{\sc A.~Aotani, M.~Mimura, and M.~Mollee}, {\em A model aided understanding of
  spot pattern formation in chemotactic e. coli colonies.}, Japan J. Ind. Appl.
  Math., 27 (2010), pp.~5--22.

\bibitem{Aranson02}
{\sc I.~S. Aranson and L.~Kramer}, {\em The world of the complex
  {G}inzburg-{L}andau equation}, Rev. Mod. Phys., 74 (2002), pp.~99--143.

\bibitem{arnold}
{\sc V.~I. Arnol{\cprime}d}, {\em Matrices depending on parameters}, Uspehi
  Mat. Nauk, 26 (1971), pp.~101--114.

\bibitem{deelanger}
{\sc E.~Ben-Jacob, H.~Brand, G.~Dee, L.~Kramer, and J.~S. Langer}, {\em Pattern
  propagation in nonlinear dissipative systems}, Phys. D, 14 (1985),
  pp.~348--364.

\bibitem{Bradley88}
{\sc R.~M. Bradley and J.~M.~E. Harper}, {\em Theory of ripple topography
  induced by ion bombardment}, Journal of Vacuum Science and Technology A:
  Vacuum, Surfaces, and Films, 6 (1988), pp.~2390--2395.

\bibitem{Bradley10}
{\sc R.~M. Bradley and P.~D. Shipman}, {\em Spontaneous pattern formation
  induced by ion bombardment of binary compounds}, Phys. Rev. Lett., 105
  (2010), p.~145501.

\bibitem{choslaut}
{\sc P.~Chossat and R.~Lauterbach}, {\em Methods in equivariant bifurcations
  and dynamical systems}, vol.~15 of Advanced Series in Nonlinear Dynamics,
  World Scientific Publishing Co. Inc., River Edge, NJ, 2000.

\bibitem{couairon}
{\sc A.~Couairon and J.-M. Chomaz}, {\em Absolute and convective instabilities,
  front velocities and global modes in nonlinear systems}, Phys. D, 108 (1997),
  pp.~236--276.

\bibitem{deelanger1}
{\sc G.~Dee and J.~S. Langer}, {\em Propagating pattern selection}, Phys. Rev.
  Lett., 50 (1983), pp.~383--386.

\bibitem{Droz00}
{\sc M.~Droz}, {\em Recent theoretical developments on the formation of
  liesegang patterns}, Journal of Statistical Physics, 101 (2000),
  pp.~509--519.

\bibitem{dumortier}
{\sc F.~Dumortier}, {\em Techniques in the theory of local bifurcations:
  blow-up, normal forms, nilpotent bifurcations, singular perturbations}, in
  Bifurcations and periodic orbits of vector fields ({M}ontreal, {PQ}, 1992),
  vol.~408 of NATO Adv. Sci. Inst. Ser. C Math. Phys. Sci., Kluwer Acad. Publ.,
  Dordrecht, 1993, pp.~19--73.

\bibitem{eckwayne}
{\sc J.-P. Eckmann and C.~E. Wayne}, {\em The nonlinear stability of front
  solutions for parabolic partial differential equations}, Comm. Math. Phys.,
  161 (1994), pp.~323--334.

\bibitem{Fenichel71}
{\sc N.~Fenichel}, {\em Persistence and smoothness of invariant manifolds for
  flows}, Indian Univ. Math. J., 21 (1972), pp.~193--226.

\bibitem{Fenichel74}
\leavevmode\vrule height 2pt depth -1.6pt width 23pt, {\em Asymptotic stability
  with rate conditions}, Indian Univ. Math. J., 23 (1974), pp.~1109--1137.

\bibitem{Fenichel77}
\leavevmode\vrule height 2pt depth -1.6pt width 23pt, {\em Asymptotic stability
  with rate conditions ii}, Indian Univ. Math. J., 26 (1977), pp.~81--93.

\bibitem{Foard12}
{\sc E.~M. Foard and A.~J. Wagner}, {\em Survey of morphologies formed in the
  wake of an enslaved phase-separation front in two dimensions}, Phys. Rev. E,
  85 (2012), p.~011501.

\bibitem{friedrich}
{\sc R.~Friedrich, G.~Radons, T.~Ditzinger, and A.~Henning}, {\em Ripple
  formation through an interface instability from moving growth and erosion
  sources}, Phys. Rev. Lett., 85 (2000), pp.~4884--4887.

\bibitem{bradleygelfand}
{\sc M.~P. Gelfand and R.~M. Bradley}, {\em Highly ordered nanoscale patterns
  produced by masked ion bombardment of a moving solid surface}, Phys. Rev. B,
  86 (2012), p.~121406.

\bibitem{Goh11}
{\sc R.~Goh, S.~Mesuro, and A.~Scheel}, {\em Spatial wavenumber selection in
  recurrent precipitation}, SIAM Journal on Applied Dynamical Systems, 10
  (2011), pp.~360--402.

\bibitem{imd}
{\sc G.~Iooss and A.~Mielke}, {\em Time-periodic {G}inzburg-{L}andau equations
  for one-dimensional patterns with large wave length}, Z. Angew. Math. Phys.,
  43 (1992), pp.~125--138.

\bibitem{katok}
{\sc A.~Katok and B.~Hasselblatt}, {\em Introduction to the modern theory of
  dynamical systems}, vol.~54 of Encyclopedia of Mathematics and its
  Applications, Cambridge University Press, Cambridge, 1995.
\newblock With a supplementary chapter by Katok and Leonardo Mendoza.

\bibitem{Keller81}
{\sc J.~B. Keller and S.~I. Rubinow}, {\em Recurrent precipitation and
  liesegang rings}, The Journal of Chemical Physics, 74 (1981), pp.~5000--5007.

\bibitem{kirchgassner}
{\sc K.~Kirchg{\"a}ssner}, {\em Wave-solutions of reversible systems and
  applications}, J. Differential Equations, 45 (1982), pp.~113--127.

\bibitem{krekhov}
{\sc A.~Krekhov}, {\em Formation of regular structures in the process of phase
  separation}, Phys. Rev. E, 79 (2009), p.~035302.

\bibitem{liese}
{\sc R.~Liesegang}, {\em {\"U}ber einige {E}igenschaften von {G}allerten},
  Naturwiss. Wochenschr., 11 (1896), pp.~353--362.

\bibitem{matsuyama}
{\sc M.~Matsushita, F.~Hiramatsu, N.~Kobayashi, T.~Ozawa, and T.~Yamazaki,
  Y.and~Matsuyama}, {\em Colony formation in bacteria: experiments and
  modeling}, Biofilms, 1 (2004), pp.~305--317.

\bibitem{Miecgl}
{\sc A.~Mielke}, {\em The {G}inzburg-{L}andau equation in its role as a
  modulation equation}, in Handbook of dynamical systems, {V}ol. 2,
  North-Holland, Amsterdam, 2002, pp.~759--834.

\bibitem{nozakibekki}
{\sc K.~Nozaki and N.~Bekki}, {\em Pattern selection and spatiotemporal
  transition to chaos in the {G}inzburg-{L}andau equation}, Phys. Rev. Lett.,
  51 (1983), pp.~2171--2174.

\bibitem{SansSchRad07}
{\sc J.~D. Rademacher, B.~Sandstede, and A.~Scheel}, {\em Computing absolute
  and essential spectra using continuation}, Physica D: Nonlinear Phenomena,
  229 (2007), pp.~166 -- 183.

\bibitem{Scheel00}
{\sc B.~Sandstede and A.~Scheel}, {\em Absolute and convective instabilities of
  waves on unbounded and large bounded domains}, Physica D: Nonlinear
  Phenomena, 145 (2000), pp.~233 -- 277.

\bibitem{gluing}
{\sc B.~Sandstede and A.~Scheel}, {\em Gluing unstable fronts and backs
  together can produce stable pulses}, Nonlinearity, 13 (2000), pp.~1465--1482.

\bibitem{SansSch01}
{\sc B.~Sandstede and A.~Scheel}, {\em On the structure of spectra of modulated
  travelling waves}, Mathematische Nachrichten, 232 (2001), pp.~39--93.

\bibitem{Scheel03}
{\sc B.~Sandstede and A.~Scheel}, {\em Defects in oscillatory media: Toward a
  classification}, SIAM Journal on Applied Dynamical Systems, 3 (2004),
  pp.~1--68.

\bibitem{standing}
{\sc B.~Sandstede and A.~Scheel}, {\em Absolute instabilities of standing
  pulses}, Nonlinearity, 18 (2005), pp.~331--378.

\bibitem{basins}
\leavevmode\vrule height 2pt depth -1.6pt width 23pt, {\em Basin boundaries and
  bifurcations near convective instabilities: a case study}, J. Differential
  Equations, 208 (2005), pp.~176--193.

\bibitem{radesherr}
{\sc M.~J. Smith, J.~D.~M. Rademacher, and J.~A. Sherratt}, {\em Absolute
  stability of wavetrains can explain spatiotemporal dynamics in
  reaction-diffusion systems of lambda-omega type}, SIAM J. Appl. Dyn. Syst., 8
  (2009), pp.~1136--1159.

\bibitem{Tobias98}
{\sc S.~Tobias, M.~Proctor, and E.~Knobloch}, {\em Convective and absolute
  instabilities of fluid flows in finite geometry}, Physica D: Nonlinear
  Phenomena, 113 (1998), pp.~43 -- 72.

\bibitem{vanSaarloos03}
{\sc W.~van Saarloos}, {\em Front propagation into unstable states}, Physics
  Reports, 386 (2003), pp.~29 -- 222.

\bibitem{vanSaarloos92}
{\sc W.~van Saarloos and P.~Hohenberg}, {\em Fronts, pulses, sources and sinks
  in generalized complex {G}inzburg-{L}andau equations}, Physica D: Nonlinear
  Phenomena, 56 (1992), pp.~303 -- 367.

\end{thebibliography}
\bibliographystyle{siam}

\end{document}